\documentclass[a4paper,11pt]{article}    
\newtheorem{theorem}{\bf Theorem}[section]
\newtheorem{lemma}[theorem]{\bf Lemma}

\newenvironment{definition}{\noindent{\bf Definition}}{\vspace{2ex}}
\newenvironment{remark}{\noindent{\bf Remark}}{\vspace{2ex}}
\usepackage{float}
\usepackage{caption}
\usepackage{graphicx}
 \usepackage{amssymb}
\begin{document}
\title{A Novel Fourier Theory on Non-linear Phases and Applications} 
\author{Tao QIAN\thanks{Department of Mathematics, FST, University of Macau, PO Box 3001, Macau, P. R. China. Supported by Macao University Multi-Year Research Grant (MYRG) MYRG2016-00053-FST and Macao Government Science and Technology Foundation FDCT 079/2016/A2. fsttq@umac.mo}}
\maketitle
\begin{abstract}{Positive time varying frequency representation for transient signals has been a hearty desire of signal analysts due to its theoretical and practical importance. During approximately the last two decades there has formulated a signal decomposition and reconstruction method rooted in harmonic and complex analysis giving rise to the desired signal representation. The method decomposes any signal into a few basic signals that possess positive instantaneous frequencies. The theory has profound relations with classical mathematics and can be generalized to signals defined in higher dimensional manifolds with vector and matrix values, and in particular, promotes rational approximation in higher dimensions. This article mainly serves as a survey. It also gives a new proof for a general convergence result, as well as a proof for the necessity of multiple selection of the parameters.

Expositorily, for a given real-valued signal $f$ one can associate it with a Hardy space function $F$ whose real part coincides with $f.$ Such function $F$ has the form $F=f+iHf,$ where $H$ stands for the Hilbert transformation of the context. We develop fast converging expansions of $F$ in orthogonal terms of the form
\begin{eqnarray*}\label{1} F=\sum_{k=1}^\infty c_k B_k,\end{eqnarray*}
where $B_k$'s are also Hardy space functions but with the additional properties
\begin{eqnarray*}\label{2} B_k(t)=\rho_k(t)e^{i\theta_k(t)},\quad \rho_k\ge 0, \quad \theta_k'(t)\ge 0, \quad {\rm a.e.}\end{eqnarray*}
The original real-valued function $f$ is accordingly expanded  \[ f=\sum_{k=1}^\infty \rho_k(t) \cos\theta_k(t)\]
which, besides the properties of $\rho_k$ and $\theta_k$ given above, also satisfies
\[ H( \rho_k\cos\theta_k)(t)= \rho_k(t) \sin\theta_k(t).\]
Real-valued functions $f(t)=\rho(t)\cos\theta(t)$ that satisfy the condition
\[ \rho\ge 0, \quad \theta'(t)\ge 0, \quad H(\rho\cos\theta)(t)= \rho(t) \sin\theta(t),\]
are called mono-components. If $f$ is a mono-component, then
the phase derivative $\theta'(t)$ is defined to be instantaneous
frequency of $f.$ The above described positive-instantaneous
frequency expansion is a generalization of the Fourier
series expansion. Mono-components are crucial to
understand the concept instantaneous frequency. We will present several most important mono-component function classes. Decompositions of signals into mono-components are called
adaptive Fourier decompositions (AFDs). We note that
some scopes of the studies on the 1D mono-components and
AFDs can be extended to vector-valued or even matrix-valued signals defined on higher dimensional
manifolds. We finally provide an
account of related studies in pure and applied mathematics,
and in signal analysis, as well as applications of the
theory found in the literature.}
\end{abstract}
\tableofcontents \pagenumbering{arabic}
\section{Introduction}
It is a common sense among analysts that  $\lq\lq$The study on the unit circle is harmonic analysis; and inside the unit circle is complex analysis" (\cite{GS}), and the same is true for a manifold and its neighborhood regions.  In general, the following mechanism may be regarded as complex analysis method of harmonic analysis. When studying analysis on the boundary of a region, say for instance, in an Euclidean space, one can at the cost of one more (or $p$-more) dimension (dimensions), imbed the region together with its boundary into the larger space, while the latter is equipped with a Cauchy complex structure, including the Cauchy Theorem, Cauchy kernel and the Cauchy formula, etc. That is, one treats the boundary of the region as a co-dimension 2 (or co-dimension $p+1$) manifold in the larger space with a Cauchy complex structure. With the complex structure one can define complex Hardy spaces consisting of $\lq\lq$good" complex holomorphic functions in the regions divided by the manifold. By $\lq\lq$good" here we in particular mean the complex-Hardy functions defined on the regions having non-tangential boundary limits as projections into the right function spaces on the manifold. Conversely, each function in an appropriate function class on the manifold can be be made to correspond with one of those Hardy-non-tangential boundary limits, the latter being called analytic signals. Those ideas appeared in the lectures of M.-T. Cheng and D.-G. Deng given in Beijing University (\cite{CD}), in the book of Gorusin translated by J.-G. Chen (\cite{Gorusin}), in works of A. McIntosh and, separately, of C. Kenig and other authors, in relation to complex Hardy spaces, singular integrals, boundary value problems and related topics on Lipschitz curves and surfaces. This article serves as a survey on the study that the author and his collaborators have been undertaking by implementing the complex analysis method to harmonic and signal analysis.

The study can be divided into two parts of which one is mono-component function theory, dealing with signals possessing a non-negative instantaneous frequency function; and the other is approximation to arbitrary function with appropriate mono-components. Note that the monomials $z^n, n=0,1,\cdots,$ are particular cases of mono-components, while the Fourier series expansion is a mono-component approximation. The study that we are going to explore is a generalization of the Fourier theory into the scope of the Beurling-Lax Theorem-related direct sum decomposition of the Hardy spaces into the forward shift and backward shift invariant subspaces.

The study at beginning was motivated by the tentative definition of the concept \emph{instantaneous frequency} (IF), or, in brief,  the \emph{frequency} function, by Gabor (\cite{Gab}).
It is still a controversial concept up to the present time.  People tend to believe that for a general signal there is a certain $\lq\lq$frequency" at each moment of time. This belief is supported by sinusoidal functions that possess constant frequencies. Justification of existence of a frequency function crucially depends on how to define IF. Unfortunately, except for the sinusoidal functions the IF concept itself appears to be paradoxical: $\lq\lq$frenquency" is the oscillation number (or in the averaging sense) per unit time duration, hence a time interval is required in order to determine it; while $\lq\lq$instantaneous" involves only a time moment.  The concept is also contradictory with the uncertainty principle. There would be no uniformly accepted theory: A great variety of engineering definitions of IF have been proposed of which, in the author's opinion, all vague and self-contradictory.  None of the theory, nor the applications, are satisfied neither by mathematicians, nor by physicians, and nor by signal analysts as (\cite{Boa,Co}). The author's view is that there may not exist an anticipated IF concept for a general signal. One can, however, propose a mathematical and conceptual definition of instantaneous frequency based on which signals can be effectively analyzed. The proposed definition of IF is based on the M\"obius transformation that, by its beauty in mathematics, has been ready to be used for the IF definition.  A coherent theory that has delicate and profound relations to classical mathematics and great potential in applications has been initialized. As a new trend of Fourier analysis it emphasizes on non-linear phase phenomenon and consists of two parts: defining the IF concept, and decomposing general signal into those possessing IF. We call the signals possessing an IF function as \emph{mono-components} (or MCs). Signals that do not possess an IF are called \emph{multi-components}. There are classical function classes belonging to the mono-component function class. There are also newly constructed interesting and significant function classes belonging to the mono-component function class.  The mono-component function theory is a combined effort by world harmonic and signal analysts (see \S 2 and the related literature in the references).

For mono-component approximation we note that G. Weiss and M. Weiss published a paper in 1962 re-proving the Nevanlinna factorization theorem in the complex Hardy spaces of one complex variable. The factorization result is a crucial tool in the complex Hardy spaces theory. Directly related to the approximation, M. Nahon, in 2000, in his Ph.D. thesis at Yale University under guidance of R. Coifman developed the non-linear phase unwinding algorithm (UWA) to expand any analytic signal into a series of Blaschke products (\cite{Na}). In their 2016 paper \cite{CS} R. Coifman and S. Steinerberger formally published the UWA theory and algorithm, and further developed some aspects initialized in \cite{Na}. A later paper by those authors together with H-T Wu developed the practical use of the unwinding method to functions of finite energy along with their computation method of the IFs (\cite{CSW}). More recently a new paper by R. Coifman and J. Peyri\'ere studies invariant subspace decompositions including the Schauder basis property of the unwinding series (\cite{CP}).  Being unaware of Nahon's thesis, Qian independently studied the UWA method and proved its $H^2$-convergence in Remark 4.4, \cite{Q} (2010), and independently uses the same English word $\lq\lq$unwinding" to call it in \cite{QLS} (2013). It is noted that UWA is a special case of UWAFD, the letter being incorporated with a sifting process using a generalized backward shifting operator together, as well as a maximal selection principle (\cite{Q,QWa}).

As already mentioned unwinding method is only one of the two main ones in the adaptive approximation methodology. The other one is the so called maximal selection principle (MSP). The terminology adaptive Fourier decomposition (AFD), as a matter of fact, at the very beginning started from the MSP type (\cite{QWa}), and further extended to the UWA type, as indeed the latter being also adaptive.
The 1D maximal selection type AFD heavily depends on the factorization properties of one complex variable. For multi-variables cases, either with the several complex variables or the Clifford algebra settings,  the AFD methods are not directly applicable. In our latest studies we extend the AFD idea to reproducing kernel Hilbert spaces with certain boundary vanishing property, called pre-orthogonal AFD, or POAFD in brief. The mono-component function theory and the related AFD approximation theory have found significant applications, including those in system identification, signal and image processing, etc. (\cite{MQ1, MQW, QZL, DQ4, ZQDM}). We will include some literature with short descriptions on engineering applications.

 The writing of the paper is organized as follows.  In \S 2 we present the main results of mono-component function theory, including the definition of mono-component function, the inner function type, the Bedrosian type, and the starlike type mono-components. In \S 3 we give an account on various kinds of AFD algorithms in the classical setting, as well as in reproducing kernel Hilbert spaces. In \S 4 we provide information on related studies and applications.
\section{Mono-component Function Theory}
\subsection{Mono-component and IF}

In 1946 Gabor proposed his \emph{analytic signal} approach (\cite{Ga}). Throughout this article we restrict ourselves to only signals with finite energy, or $L^2$-functions. The theory on the unit circle is parallel with that on the real line. To explain the idea we most time restrict ourselves to the unit circle case. We occasionally jump into the upper-half space context, including for instance when we describe the ideas in relation to the Bedrosian type results in terms of Fourier transform.  Let $s(t)$ be a real-valued signal of finite energy on the unit circle $\partial {\bf D}$ where ${\bf D}$ denotes the unit disc. The associated analytic signal, denoted by $s^+(t),$ is defined
\begin{eqnarray}\label{s+} s^+(e^{it})=\frac{1}{2}\left( s(e^{it})+i \tilde{H}s (e^{it})+c_0\right),\end{eqnarray}
where $\tilde{H}$ is the circular Hilbert transformation, and $c_0$ is the $0$-th Fourier coefficient, or average of $s$ on the circle. That is
\[\tilde{H}s(e^{it})=\frac{1}{2\pi}{\rm p. v.}\int_{0}^{2\pi}f(e^{iu})\cot \left(\frac{t-u}{2}\right)du, \quad c_0=\frac{1}{2\pi}\int_0^{2\pi}f(e^{iu})du.\]
We note  $s^+$ is the non-tangential boundary limit of the Cauchy integral of $s$ (the Plemelj formula):
\[ s^+(e^{it})=\lim_{z\to e^{it}}\frac{1}{2\pi}\int_0^{2\pi}\frac{f(e^{iu})}{z-e^{iu}}e^{iu}du, \qquad {a.e.}\]
  The fact that $s$ is real-valued makes the Hilbert transform part $\frac{1}{2}\tilde{H}f$ in (\ref{s+}) the purely imaginary part of $s^+,$ and $s=2{\rm Re}s^+-c_0.$
There also holds the following relation that in the real line context corresponds to the Laplace transform
\[ s^+(e^{it})=\sum_{k=0}^\infty c_ke^{ikt}.\]
What is important is that $s^+(e^{it})$ has a holomorphic continuation into the interior of the disc
\[  s^+(z)=\sum_{k=0}^\infty c_kz^{k}, \quad |z|<1, \] as a Hardy $H^2({\bf D})$ function in the sense that whose non-tangential boundary limit coincides with $s^+(e^{it}).$
The Fourier multiplier of the circular Hilbert transformation is $-i{\rm sgn},$ that is, if
\[ s(e^{it})=\sum_{k=-\infty}^\infty c_ke^{ikt}, \quad c_k=\frac{1}{2\pi}\int_0^{2\pi}s(e^{iu})e^{-iku}du,\]
where ${\rm sgn}(k)=1,$ if $k>0,$ and ${\rm sgn}(k)=-1,$ if $k<0,$ and ${\rm sgn}(0)=0,$ then
\[ \tilde{H}s(e^{it})=\sum_{k=-\infty}^\infty (-i){\rm sgn}(k)c_ke^{ikt}.\]
In the sequel we drop the tilde sign above $\tilde{H}$ and write it simply as $H.$

This Fourier multiplier form of the Hilbert transform gives rise to the Hilbert transform characterization of the Hardy spaces. If restricted to the $L^2$ cases, it is: A function $s$ of finite energy belongs to the Hardy $H^2$ space if and only if $Hs=-is$ (\cite{Q17}).  This result holds in general contexts including the upper-half space cases in one and higher dimensions (\cite{DMQ1,DLQ1,DLQ2}).

In writing $s^+(e^{it})=\rho (t)e^{i\theta (t)},$ Gabor defined that the derivative of the phase function, $\theta'(t),$ to be the instantaneous frequency of $s(e^{it}).$ In commenting on this definition we would say that the definition is $\lq\lq$good", because if we take the example that for a positive integer $n, s(e^{it})=\cos (nt),$ then in such way, $s^+(e^{it})=e^{int},$ and the phase derivative is $n$, being complementary with the common sense. Gabor's definition, however, is not valid for general signals $s\in L^2(\partial {\bf D}),$ but only tentative, due to the following reasons. First $s,$ and thus $s^+$ as well, is an equivalent class of Lebesgue square-integrable functions that cannot be expected to be smooth and thus have phase derivative; and secondly, the derivative, if exists, cannot be expected to be non-negative, as required in physics, and thus cannot stand as a qualified instantaneous frequency function. It is the signal analysts who decided that the IFs should be non-negative and thus can be effectively analyzed with applications.  The primary importance is that the instantaneous frequency concept is generated from physics practice: it is an extension of the vibrating frequency. In the average sense the phase derivative of an analytic signal is non-negative as read out from the relation
\[ \frac{1}{2\pi}\int_{0}^{2\pi}\theta'(t)|s^+(e^{it})|^2dt=\sum_{k=0}^\infty k|c_k|^2,\] (\cite{Co,DQY}).
Pointwisely, however, the phase derivative of an analytic signal  can be negative. For instance, for any non-trivial outer function in the complex Hardy space we have a set of positive Lebesgue measure on which the phase derivatives are strictly less than zero (\cite{Q20}).

The strategy is to define a function set that exactly contains the signals having well defined non-negative analytic phase derivatives. It is the set of the functions called mono-components.

\begin{definition} \label{mono}(Tao Qian 2006 \cite{Q19}) Let $s$ be a real- or complex-valued signal on the unit circle of finite energy. We call $s$ a mono-component, or real-mono-component, if its analytic signal, or equivalently its projection into the Hardy space $H^2,$ viz., $s^+(t)=\frac{1}{2}\left( s(t)+i Hs (t)+c_0\right),$ in its phase-amplitude representation $s^+(t)=\rho (t)e^{i\theta (t)}$ satisfies $\theta'(t)\ge 0,$ a.e., where the phase derivative $\theta'(t)$ is defined through the non-tangential limit of the same quantity from inside of the region. Precisely, in the unit circle case,
\[ \theta'(t)=\lim_{r\to 1-} \theta'_r(t), \quad {\rm a.e.},\]
where $s^+(re^{it})=\rho_r(t)e^{i\theta_r(t)}$ is the holomorphic continuation inside the unit disc.
When $s$ is a mono-component we call $s^+$ a complex-mono-component, or simply mono-component as well.
When and only when $s$ is a mono-component it has an instantaneous frequency function defined as its non-negative analytic phase derivative $\theta'(t).$ \end{definition}

Since $s^+$ is the non-tangential boundary limit of a Hardy space function inside the unit disc, $\theta'_r(t)$ everywhere exists, and
\[ \theta'_r(t)={\rm Re}\{\frac{re^{it}{s^+}'(re^{it})}{s^+(re^{it})}\}.\]
We note that the class of mono-component functions is closed under the multiplication operation but not the addition.

\subsection{The Inner Function Type Mono-components}

The above definition allows a number of interesting mono-component subclasses. First we will mention the class of inner functions. It is easily observed that the boundary function of the canonical M\"obius transform mapping $a\in {\bf D}$ to zero, that is
\[ e_a(e^{it})=\frac{e^{it}-a}{1-\overline{a}z}=e^{i\theta_a(t)},\] is an analytic signal whose phase derivative $\theta'_a$ is the Poisson kernel for the disc (\cite{Ga,QCL,Q18}). This early study along this direction was joined by Qiu-Hui Chen and Luo-Qing Li. This implies that finite Blaschke products (Blaschke products with finitely many zeros) are all mono-components. The question is whether infinite Blaschke products are mono-components. As an application of the Julia-Wolff-Carath\'eodory Theorem the following result on general inner functions (containing finite and infinite Blaschke products and singular inner functions) is proved (\cite{Q20}).

\begin{theorem} (Tao Qian 2009)\label{inner} Let $\theta$ be a real-valued Lebesgue measurable function on the unit circle. Then the phase function $e^{i\theta}$ is a mono-component if and only if $e^{i\theta}$ is the non-tangential boundary limit of an inner function, or, equivalently, if and only if $H(e^{i\theta})=-ie^{i\theta}.$ \end{theorem}

The earlier study in \cite{Daniel} gives good observations and partial results. It is noted that in the earlier digital signal processing (DSP) literature the fact that Blaschke products possess positive phase derivative functions were stated without valid proof, as far as being aware by the author (\cite{ChengQS}). DSP scholars and engineers have been using the concept physically realizable signals with minimum phase that are outer functions without rigorous proof either. The minimum phase phenomenon is based on the fact that all boundary values of inner functions have non-negative phase derivative. The reference \cite{Q20} proves the opposite property for outer functions: Under mild conditions that guarantee existence of the phase derivative function $\theta'(t)$ for an outer function there holds
\[ \int_0^{2\pi} \theta'(t)dt=0.\]

The inner and outer functions are thus characterized in terms of sign properties of their phase derivatives.

\subsection{The Bedrosian Type Mono-components} The second class of mono-components is called the Bedrosian type. The
Bedrosian theorem asserts the relation
\[ H(fg)=fHg\] under two forms of conditions, both being based on Fourier spectrum property of the functions. The first form of the conditions is that there exists $\sigma>0$ such that ${\rm supp}\hat{f}\subset [-\sigma,\sigma],$ and ${\rm supp}\ \hat{g}\subset (-\infty,\sigma]\cap [\sigma,\infty).$ The second form of the conditions is that both functions $f$ and $g$ are in the Hardy $H^2$ space. In the language of the Fourier spectrum, by recalling the Paley-Wiener Theorem for the Hardy space functions the second form of the conditions is equivalent with $f,g\in L^2$, ${\rm supp} \hat{f}\subset [0,\infty)$ and  ${\rm supp}\ \hat{g}\subset [0,\infty)$. The idea of using the Bedrosian type results is as follows. Suppose that $e^{i\theta}$ is an analytic signal with the property $\theta'(t)\ge 0,$ a.e. This kind of functions now have all been characterized by Theorem \ref{inner}.  One wishes to find a non-negative function $\rho(t)$ that makes the Bedrosian type relation  $H(\rho(t)e^{i\theta}=\rho H(e^{i\theta})$ hold. For such a function $\rho$  there holds
\[ H(\rho e^{i\theta})=\rho H(e^{i\theta})=(-i)\rho e^{i\theta}.\]
By recalling the Hilbert transform characterization of the Hardy space functions the last equality implies that $\rho e^{i\theta}$ is an analytic signal, and, due to the positivity of the phase derivative $\theta'(t),$ it is a mono-component.

However, the original Bedrosian theorem cannot be directly used. The first form of the conditions refers to bandlimiting properties of the functions $f$ and $g.$ That, unfortunately, are not our case: The inner functions $g$ do have the full spectrum range. The second form of the conditions would require that the amplitude function itself is the boundary limit of some Hardy space function. But it is not the case either (see the example given in (\ref{exa})).

In order to enrich the mono-component class new conditions for the Bedrosian identity to hold were seeking by mainly a group of Chinese harmonic and signal analysts, including Qiu-Hui Chen, Li-Hui Tan, Rui Wang, Si-Lei Wang, Yue-Sheng Xu, Dun-Yan Yan, Li-Hua Yang, Bo Yu, Hai-Zhang Zhang, Li-Xin Yan, etc., as well as the author (\cite{Xu,YuZ,QWYZ,Tan-Shen-Yang,WangSL,QT1}, etc.), with Fourier analysis methods and complex analysis methods. The most successful result along this line is based on the following observation.

The essential structure of Bedrosian type mono-components is as follows:
\begin{eqnarray}\label{exa} s(e^{it})=\left(\frac{1}{1-\overline{a}_1 e^{it}}+\frac{1}{1-a_1e^{-it}}\right)\frac{e^{it}-a_1}{1-\overline{a}_1e^{it}}
\frac{e^{it}-a_2}{1-\overline{a}_2e^{it}}.\end{eqnarray}
On the circle it is a real-valued function multiplied with an order-2 Blaschke product.
In verifying that $s(e^{it})$ is a Bedrosian type mono-component, the key point is that
\[ \frac{1}{1-\overline{a}_1 z}\frac{z-a_1}{1-\overline{a}_1z}
\frac{z-a_2}{1-\overline{a}_2z}\]
is an analytic function in the disc; and for $|z| = 1,$ the product
\[\frac{1}{1-{a}_1\overline{z}}\frac{z-a_1}{1-\overline{a}_1z}
\frac{z-a_2}{1-\overline{a}_2z}\]
has an analytic continuation to the interior part of the disc. As result, s(z) is
a bounded analytic function. Since $\frac{1}{1-\overline{a}_1 e^{it}}+\frac{1}{1-a_1e^{-it}}$ is real-valued
and has finitely many sign-change pints on $|z| = 1,$ it is, therefore, a so called \emph{generalized amplitude} on the circle.
We have the following general result (\cite{QT1}, the finite order Blaschke products case is proved in \cite{TYH}).

\begin{theorem}\label{Bedrosian} Let $\phi (e^{it})$ be an infinite Blaschke product, where $a_1,\cdots,a_n,\cdots$ are the totality of its zeros, the multiples being all counted. Then (1) $\rho (t)$ is a real-valued function such that $\rho (t)\phi (e^{it})\in H^p(\partial {\bf D}), 1\leq p\leq \infty,$ if and only if $\rho$ is the real part of some function in the backward-shift-invariant subspace induced by the Blaschke product $\phi (e^{it}),$ that is $\rho\in {\rm Re}\{H^p(\partial {\bf D})\cap \phi (e^{it})\overline{H^p(\partial {\bf D})}\};$
and (2) For $1<p<\infty,$ $\rho\in {\rm Re}\{H^p(\partial {\bf D})\cap \phi (e^{it})\overline{H^p(\partial {\bf D})}\}$ if and only if, in the $L^p$ norm sense,
\[ \rho(t)={\rm Re}\{\sum_{k=1}^\infty c_kB_k(e^{it})\},\]
where $c_k=\langle \rho (t), B_k(e^{it})\rangle=\frac{1}{2\pi}\int_0^{2\pi}\rho (t)\overline{B_k(e^{it})}dt, k=1,2,\cdots,$
and $\{B_k\}_{k=1}^\infty$ is the rational orthonormal system (or TM-system) generated by $a_1,\cdots,a_k,\cdots,$ that is
\begin{eqnarray}\label{TM} B_k(z)=\frac{\sqrt{1-|a_k|^2}}{1-\overline{a}_kz}\prod_{l=1}^{k-1}\frac{z-a_l}{1-\overline{a}_lz}.\end{eqnarray}
\end{theorem}

\subsection{The Non-Bedrosian Type Mono-components: The Starlike and Boundary Starlike Type}

The third type of mono-components is called the non-Bedrosian type which contains all $p$-starlike as well as boundary starlike functions in one complex variable. This kind of mono-components explores a different type of connections between mono-components and conformal mappings.  Let $f$ denote a univalent conformal mapping that, with $f(0)=0,$ maps the unit disc together with its continuous and thus rectifiable boundary. Obviously $f(e^{it})$ is a complex mono-component as its phase function is increasing along with increasing of the angular variable $t.$ Below we will denote by ${\cal S}^*$ the set of such starlike functions. Next we define several other likewise function classes including $p$-starlike functions as follows.

\begin{definition} Let $p$ be any positive integer. Denote by ${\cal S}(p)$ the set of $p$-valent holomorphic functions satisfying the following conditions:\\
 (i) There exists $r: 0<r<1,$ such that for all $z: r<|z|<1,$ there holds
 ${\rm Re}\{\frac{zf'(z)}{f(z)}\}>0;$ and\\
 (ii) $\int_0^{2\pi}{\rm Re}\{\frac{zf'(z)}{f(z)}\}dt=2p\pi$ for all $z: r<|z|<1.$ Functions belonging to ${\cal S}(p)$ are called $p$-starlike functions.
 \end{definition}

 \begin{definition} A function $f$ is said to be a weak $p$-valent starlike function, and denoted $f\in {\cal S}_w(p)$, if and only if it is holomorphic in ${\bf D}$ with precisely $p$ zeros in ${\bf D}$ (including multiples) and with the expression
 \[ f(z)=[h(z)]^p\prod_{k=1}^p\frac{(z-a_k)(1-\overline{a}_kz)}{z},\]
 where $h\in {\cal S}^*.$\end{definition}

With $p=1$ and $a_1=0$ we obtain ${\cal S}_w(1)={\cal S}(1).$ The article \cite{Hummel} shows that ${\cal S}(p)$ is a proper subset of ${\cal S}_w(p)$. The advantage of the latter is that functions in ${\cal S}_w(p)$ have an explicit representation formula. In order to reveal the essential structure we assume the convenient property that functions under study have a holomorphic continuation to an open neighborhood of the closed unit disc. Denote by ${\cal A}$ the set of such holomorphic functions, one can show ${\cal A}\cap {\cal S}(p)={\cal A}\cap {\cal S}_w(p)$ (\cite{QT1}). To describe the relation between mono-components and the starlike function family we need two more definitions.

\begin{definition} (\cite{Ly}) A univalent function is said to be a boundary starlike function with respect to the origin if $f$ is holomorphic in ${\bf D}, \ \lim_{r\to 1-}f(r)=0, \ f({\bf D})$ is starlike with respect to the origin, and ${\rm Re}\{e^{i\alpha}f(z)\}>0$ for some real number $\alpha$ and all $z\in {\bf D}.$ Denote by ${\cal G}^*$ the set of all boundary starlike functions with respect to the origin. \end{definition}

The following definition specifies a class of mono-components.

\begin{definition}
Let $f(e^{it})=\rho (t)e^{i\theta (t)}\in L^p(\partial {\bf D}), p\ge 1.$ Then $f$ is called a Hilbert-{n}, or $H$-$n$ atom, if it satisfies the following conditions:\\
(1) $H(\rho\cos\theta)=\rho\sin\theta;$\\
(2) $\rho\ge 0, \theta'\ge 0$ a.e.; and\\
(3) $\int_0^{2\pi} \theta'(t)dt=n\pi.$
\end{definition}

Note that due to (1) $f$ has a holomorphic continuation into the unit disc as a Hardy space function. In (2) the phase derivative $\theta'$ takes the sense given in \ref{mono}. The condition (3) refers to the multivalent degree of $f.$ The concept $H$-p atom was first proposed in \cite{Q19} for $p=2$ with the result that a function $f$ is a $H$-2 atom if and only if $f$ is a starlike function about the origin. Some further studies along this line for $p=2n$ are give in \cite{TYH}. The following result ultimately presents the relation between the $H$ atoms and the starlike-boundary starlike functions.

\begin{theorem}
Assume that $f$ is holomorphic in $\overline{\bf D}$ having $p$ zeros in the open disc ${\bf D}.$ Then $f(e^{it})$ is a $H$-$n$ atom, $n\ge 1,$ if and only if
\[ f^2(z)=\left[\prod_{i=1}^p h_i(z)\right]^2\prod_{j=1}^{n-2p}g_j^2(z)=\left[\prod_{k=1}^p(z-a_k)\left(\frac{1}{z}-\overline{a}_k\right)\right]^2
\left[\prod_{k=1}^{n-2p}(z-b_k)\left(\frac{1}{z}-\overline{b}_k\right)\right][h(z)]^n,\]
where $\{a_k\}_{k=1}^p$ are the zeros of $f(z)$ inside the unit disc, $\{b_k\}_{k=1}^{n-2p}$ are the zeros of $f(z)$ on the unit circle (both can be with multiples), $h(z)\in {\cal S}^*, h_i\in {\cal S}_w(1),$ and $g_j(b_jz)\in {\cal G}^*$ are all holomorphic in $\overline{\bf D},\ i=1,\cdots,p,\ j=1,\cdots,n-2p.$
\end{theorem}

The results on mono-component functions in all the three categories, viz., the inner function type, the Bedrosian type, and the starlike type, are not only important results in the mono-component function theory, but also new understanding to related topics in the classical harmonic and complex analysis.


\section{Adaptive Fourier Approximations}

In this part we will give expository descriptions of adaptive Fourier approximations (AFD). In the one complex variable cases AFD gives rise to positive frequency expansions of signals into holomorphic rational functions, while in higher dimensions AFD offers at the moment at least fast converging rational approximations.

\subsection{Mono-component Decomposition of Signals in General}

The idea of positive frequency decompositions of signals is not new, it goes back to more than two hundreds years ago in relation to the name Jean Baptiste Joseph Fourier and other names. Fourier series will be exactly a particular case of the general theory that we are now to present.

Let $s$ be a real-valued function defined on the unit circle $\partial {\bf D}$ with finite energy. We recall that its Hardy $H^2$ space projection is $s^+=\frac{1}{2}(s+iHs+c_0).$ The simple relation
$s=2{\rm Re}\{s^+\}-c_0$ implies that a complex-mono-component decomposition $s^+(e^{it})=\sum_{k=1}^\infty \rho_k(t)e^{i\theta_k(t)} $ (The possible coefficients $c_k$ are absorbed by the corresponding amplitudes $\rho_k$ and the phases $\theta_k$) will gives rise to a real-mono-component decomposition, or positive-frequency decomposition of $s:$ $s(e^{it})=-c_0+\sum_{k=1}^\infty \rho_k(t)\cos\theta_k(t).$ We are hence reduced to decomposing the complex Hardy space function $s^+.$

The philosophy  is to find intrinsic constructing blocks of positive-time varying-instantaneous frequency. Here $\lq\lq$intrinsic" has profound meaning. Under the present setting we understand it as $\lq\lq$fast converging". In \cite{QHLW} we show that for any Hardy space function $s^+$ and any $\epsilon >0,$ there exist a constant $c$ and two $1$-starlike functions $m_1$ and $m_2$ such that
\[ \|s^+-(c+m_1+m_2)\|\leq \epsilon.\]
The two starlike functions $m_1$ and $m_2$ are not unique, and, according to their construction, are very irregular. Deeper and more thorough consideration suggests that one would better use function systems of a certain type consisting of well behaved mono-components with explicit expressions such as rational functions, and better with monotonously increasing frequencies. In order to obtain intrinsic decomposition fast convergence would be desired.

\subsection{One Dimensional Core-Adaptive Fourier Decomposition (Core-AFD) and its Variations}

It did not take much time for the researcher to decide to use the rational orthonormal system or TM system. The difference with the traditional use of the system is that now it has to be adaptive: One must select parameters for the intrinsic construction concern, that means fast convergence. The system was already been introduced in Theorem \ref{Bedrosian}. We note that TM systems in general cannot be avoided for they are Gram-Schmidt (G-S) orthogonalization of the partial fractions with poles outside the closed unit disc, the latter being fundamental constructive building blocks of rational functions in the Hardy spaces. TM systems consist of functions of positive frequency due to their construction in (finite) Blaschke products.

In the sequel we change our function notation $s^+$ in the Hardy $H^2({\bf D})$ to $f.$ In the unit circle context  we have
 $f(z)=\sum_{l=1}^\infty c_lz^l, \sum_{l=1}^\infty |c_l|^2<\infty.$ Now we seek a decomposition of $f$ into a TM system with adaptively selected parameters.
  The collection of the functions
 $$e_a(z)=\frac{\sqrt{1-|a|^2}}{1-\overline{a}z}, \qquad a\in {\bf D},$$ consists of
 normalized Szeg\"o kernels of the disc.  Set $f=f_1.$ First write
 \[ f(z)=\langle f_1, e_{a_1}\rangle e_{a_1}(z)+\frac{f_1(z)-\langle f_1, e_{a_1}\rangle e_{a_1}(z)}{\frac{z-a_1}{1-\overline{a}_1z}}\frac{z-a_1}{1-\overline{a}_1z}.\]
 We note that in this stage $a_1$ can be any complex number in the unit disc and the above is an identity.
 Denoting
\[f_2(z)=\frac{f_1(z)-\langle f_1, e_{a_1}\rangle e_{a_1}(z)}{\frac{z-a_1}{1-\overline{a}_1z}},\] calling it the {\it reduced remainder,} the identity is re-written as
 \begin{eqnarray}\label{maximum sifting}  f(z)=\langle f_1, e_{a_1}\rangle e_{a_1}(z)+f_2(z)\frac{z-a_1}{1-\overline{a}_1z},\end{eqnarray}

We call the operator mapping $f_1$ to $f_2$ the \emph{generalized $a_1$-backward shift operator} and $f_2$ \emph{the generalized $a_1$-backward shift of} $f_1.$  The terminology is a generalization of the classical backward shift operator
\[ S(f)(z)=a_1+a_2z+\cdots +c_{k+1}z^k+\cdots  =  \frac{f(z)-f(0)}{z}.\]
Recognizing that $f(0)=\langle f, e_0\rangle e_0(z),$ the operator $S$ is the generalized $0$-backward shift operator.

Notice that the Szeg\"o kernel is the Cauchy kernel under the arc-length measure and thus has the reproducing kernel property. Due to the orthogonality in the Hilbert space and the modular one property of the M\"obius transform, we have the energy relations
\[ \|f\|^2=\| \langle f_1, e_{a_1}\rangle e_{a_1}\|^2+\|f_2\|^2 =(1-|a_1|^2)|f_1(a_1)|^2 + \|f_2\|^2.\]

 The purpose now is to extract the maximal energy portion from the term  $\langle f_1, e_{a_1}\rangle e_{a_1}(z).$  It is reduced to maximize $(1-|a_1|^2)|f_1(a_1)|^2$ among all $a_1\in {\bf D}.$ Although ${\bf D}$ is an open set we are fortunate that there exists $a_1$ in ${\bf D}$ such that
\[ a_1=\arg \max \{ (1-|a|^2)|f_1(a)|^2\ : \ a\in {\bf D}\}\]
(\cite{QWa}). The existence of such maximal selection is called \emph{Maximal Selection Principle}.  Selecting such $a_1$ and repeating the process for $f_2,$ and so on. We call the process from $f_1$ to get $f_2$ through a maximal selection $a_1$ as {\it maximal sifting from} $f_1$ \emph{to} $f_2$  \emph{through} $a_1.$ After $n$ siftings one gets
\[ f(z)=\sum_{k=1}^n \langle f_k, e_{a_k}\rangle B_k(z) + f_{n+1}\prod_{k=1}^n\frac{z-a_k}{1-\overline{a}_kz},\]
where for $k=1,\cdots ,n,$
\[ a_k=\arg \max \{ (1-|a|^2)|f_k(a)|^2\ : \ a\in {\bf D}\},\] \[B_k(z)=B_{\{a_1,\cdots,a_k\}}(z)=\frac{\sqrt{1-|a_k|^2}}{1-\overline{a}_kz}\prod_{l=1}^{k-1}
\frac{z-a_l}{1-\overline{a}_lz},\]
and, for $k=2,...,n+1,$ $f_k$ is the maximal sifting of $f_{k-1}$ through $a_{k-1},$ that is,
\[f_k(z)=\frac{f_{k-1}(z)-\langle f_{k-1}, e_{a_{k-1}}\rangle e_{a_{k-1}}(z)}{\frac{z-a_{k-1}}{1-\overline{a}_{k-1}z}}.\]
 We have the following convergent theorem.
\begin{theorem}\label{AFD}
\emph{For any  give function $f$ in the Hardy $H^2$ space, by applying the maximum sifting process at each step we have}
\[  f(z)=\sum_{k=1}^\infty \langle f_k, e_{a_k}\rangle B_k(z).\]\end{theorem}
This result was first proved in \cite{QWa} based on the complex modula $1$ property of the M\"obius transform. Below we provide a new proof releasing the modular $1$ requirement for the system functions but only based on maximal selections of the parameters. The essence of the proof is contained in several proofs of \cite{Qbook}, also see \cite{QD}.

\noindent{\bf Proof}. We prove the convergence by contradiction. Assume that through a sequence of maximally selected parameters ${\bf a}=\{a_1,\cdots,a_n,\cdots\}$ we arrive
\begin{eqnarray}\label{general} f=\sum_{k=1}^\infty \langle f_k,e_{a_k}\rangle B_k + h, \qquad h\ne 0.\end{eqnarray}

The routine argument by using the Riesz-Fisher Theorem shows that both the functions $\sum_{k=1}^\infty \langle f_k,e_{a_k}\rangle B_k $ and $h$ are in $H^2.$ We note that from the Hilbert space property $h$ is orthogonal with all $B_k,$ as well as with $\sum_{k=1}^\infty \langle f_k,e_{a_k}\rangle B_k.$

The relation (\ref{general}) can be re-written
\[ f=\left(\sum_{k=1}^M + \sum_{k=M+1}^\infty\right) \langle f_k,e_{a_k}\rangle B_k + h,\]
where by our notation,
\[ g_{M+1}=\sum_{k=M+1}^\infty \langle f_k,e_{a_k}\rangle B_k + h =G_{M+1}+h.\]
To proceed we note that \begin{eqnarray}\label{three}\langle f_k,e_{a_k} \rangle=\langle f,B_k \rangle=\langle g_k,B_k \rangle,\end{eqnarray} where $g_k=f-\sum_{l=1}^{k-1}\langle f_l,e_{a_l} \rangle B_l$ is the $k$-th {\it standard remainder}.

Therefore, we have
\[ g_{M+1}=\sum_{k=M+1}^\infty \langle g_k,B_k\rangle B_k + h.\]
Due to the density of the function set $\{e_a\}_{a\in {\bf D}}$ in $H^2,$ there exists $a\in {\bf D}$ such that $\delta \triangleq |\langle h,e_a\rangle |>0.$ We can in particular choose $a$ to be distinguished from all the selected $a_k$'s.  We are now to explore a contradiction in relation to the selections of $a_{M+1}$ for large $M.$
Now, on one hand, by the Bessel inequality applied to the infinite series part in (\ref{general}), we have
\begin{eqnarray}\label{back} |\langle g_{M+1},B_{M+1}\rangle|\to 0, \qquad {\rm as}\quad M\to 0.\end{eqnarray}

On the other hand we will show, for large $M$,
\begin{eqnarray} \label{backback}|\langle g_{M+1},B_{M+1}^{a}\rangle|> \frac{\delta}{2}.\end{eqnarray}
This is then clearly a contradiction.

The rest part of the proof is devoted to showing (\ref{backback}). Due to the relations
\begin{eqnarray}\label{g+h} |\langle g_{M+1}, B_{M+1}^a\rangle|\ge |\langle h, B_{M+1}^a\rangle|-|\langle G_{M+1}, B_{M+1}^a\rangle|\end{eqnarray} and
\begin{eqnarray}\label{h} |\langle G_{M+1}, B_{M+1}^a\rangle|\leq \|G_{M+1}\|\to 0, \qquad {\rm as}\quad M\to \infty,\end{eqnarray}
for large $M$ the lower bounds of $|\langle g_{M+1}, B_{M+1}^a\rangle|$ depend on the quantity of $|\langle h, B_{M+1}^a\rangle|$.
Now for any positive integer $M$ denote $X_{M+1}^{a}$ the $(M+1)$-dimensional space spanned by $\{e_{a}, B_1,\cdots,B_M\}.$ We have two ways to compute the energy of the projection of $h$ into  $X_{M+1}^{a},$ denoted $\|h/X_{M+1}^a\|^2.$ One way is based on the orthonormality of $\{B_1,\cdots,B_M,B_{M+1}^a\}.$ In such way, due to the orthogonality of $h$ with $B_1,\cdots,B_M,$ we have
\[ \|h/X_{M+1}^a\|^2=|\langle h,B_{M+1}^a\rangle|^2.\]
The other way is based on the orthonormalization of the $(M+1)$-tuple in the order $\{e_a,B_1,\cdots,B_M\}.$ Then we have
\[ \|h/X_{M+1}^a\|^2\ge |\langle h,e_a\rangle|^2=\delta^2.\]
Hence we have, for any $M$, $|\langle h,B_{M+1}^a\rangle|\ge \delta.$ In view of this last estimation and (\ref{h}),\ (\ref{g+h}), we arrive at the contradiction given by (\ref{backback}),\ (\ref{back}). The proof is complete.

\begin{remark}  It is noted that the selected parameters $a_1,\cdots,a_n,\cdots$ according the maximal principle may not satisfy the hyperbolic non-separable condition
 \[ \sum_{k=1}^\infty (1-|a_k|)=\infty\]
 and thus the generated TM system $\{B_k\}$ may not be a basis. By doing such decomposition one is not interested in whether the resulted TM system is a basis, but only in whether it can effectively expand the given signal $f.$  One is indeed able to do so, and, in fact, achieves fast convergence.\end{remark}

\begin{remark} For arbitrary selections of $a_1,...,a_n,...,$ we arrive a pre-mono-component decomposition: All entries in the infinite sum after being multiplied by $e^{it}$ become mono-components. If we choose $a_1=0,$ then all $B_k$'s are mono-components, and AFD offers a mono-component decomposition. \end{remark}

\begin{remark}  AFD is different from any existing greedy algorithm (\cite{DMA}) for the following reasons. (i) On the algorithm side AFD is incorporated with a generalized backward-shift operation, or a sifting process, that changes the standard remainder to reduced remainder; (ii) The reduced remainder allows multiple selections of parameters for optimal approximation at each of the iteration steps; and (iii) The backward-shift operation automatically generates an orthonormal system, or a TM system without using the Gram-Schmidt process. \end{remark}

\begin{remark} Restricted to a practical subclass the convergence rate for AFD is $M/\sqrt{n}$ where $n$ is the order of the AFD partial sum. One has to note that this is a good convergence rate for it is for non-smooth functions (boundary limits of Hardy space functions). The classical convergence theorems may look better but they are for smooth functions. \end{remark}

\subsection{Unwinding AFD (UWAFD)}

 Let $f=hg,$ where $f,g$ are Hardy $H^2({\bf D})$ functions, and $h$ is an inner function. Let $f$ and $g$ be expanded into their respective Fourier series, viz.,
\[ f(z)=\sum_{k=0}^\infty c_kz^k, \qquad g(z)=\sum_{k=0}^\infty d_kz^k.\]
The Plancherel Theorem and the
modula $1$ property of inner functions assert that
\[ \sum_{k=0}^\infty |c_k|^2=\| f\|^2=\| g\|^2=\sum_{k=0}^\infty |d_k|^2.\]
In digital signal processing (DSP) there is following result: for any $n,$
\[ \sum_{k=n}^\infty |c_k|^2 \ge \sum_{k=n}^\infty |d_k|^2\]
(see, for instance \cite{DQ2, ChengQS}).\\

In DSP this is referred to be energy front loading property of minimum phase signals among physically realizable signals. This amounts to saying that after factorizing out the inner function factor,  the convergence rate of the Fourier series of the remaining outer function speeds up. This suggests that the AFD process would be better to incorporate a factorization process for speeding up the convergence. This is reasonable: when a signal by its nature is of high frequency, one should first perform $\lq\lq$unwending" before extracting out from it a maximal portion of lower frequency. We proceed it as follows (\cite{Q}, \cite{QLS}). First we do factorization $f=f_1=I_1O_1,$ where $I_1$ and $O_1$ are, respectively, the inner and outer factors of $f.$ The factorization is based on Nevanlinna's factorization theorem, also see \cite{WW}. The outer function has the explicit integral representation
\[ O_1(z)=e^{\frac{1}{2\pi}\int_0^{2\pi}\frac{e^{it}+z}{e^{it}-z}\log |f_1(e^{it})|dt}.\]
The boundary value of the outer function is computed by using the boundary value of $f_1.$  On the boundary the above integral is taken to be of the principal integral sense. The imaginary part of the integral reduces to the circular Hilbert transform of $\log |f_1(e^{it})|.$ Next, we do a maximum sifting to $O_1.$ That gives
\[ f(z)=I_1(z)[\langle O_1, e_{a_1}\rangle e_{a_1}(z)+f_2(z)\frac{z-a_1}{1-\overline{a}_1z}],\]
where $f_2$ is the maximal shifting of $O_1$ through $a_1:$
\[ f_2(z)=\frac{O_1(z)-\langle O_1, e_{a_1}\rangle e_{a_1}(z)}{\frac{z-a_1}{1-\overline{a}_1z}}.\]
By factorizing $f_2$ into its its inner and outer factors, $f_2=I_1O_2,$ we have
\[  f(z)=I_1(z)[\langle O_1, e_{a_1}\rangle e_{a_1}(z)+I_2(z)O_2(z)\frac{z-a_1}{1-\overline{a}_1z}].\]
We next proceed a maximum sifting to $O_2,$ and so on. In such way we arrive at the unwinding AFD decomposition (\cite{Q}):

\begin{theorem}\label{unwending}
The above procedure gives rise to the unwinding AFD (UWAFD) decomposition
\[ f(z)=\sum_{k=1}^n \prod_{l=1}^k I_l(z) \langle O_k,e_{a_k}\rangle B_k(z) + f_{n+1}(z)\prod_{k=1}^n\frac{z-a_k}{1-\overline{a}_kz}\prod_{l=1}^{n}I_l(z),\]
where $f_{k+1}=I_{k+1}O_{k+1}$ is the maximal shifting of $O_{k}$ through $a_k, k=1,...,n,$ and
$I_{k+1}$ and $O_{k+1}$ are respectively the inner and outer functions of $f_{k+1}.$ Furthermore,
\[ f(z)=\sum_{k=1}^\infty \prod_{l=1}^k I_l(z) \langle O_k,e_{a_k}\rangle B_k(z).\]\end{theorem}

\begin{remark} Like AFD,  unwending AFD (UWAFD) is a mono-component or pre-mono-component decomposition. Experiments show that among various AFD type algorithms UWAFD seems to converge most rapidly, that is in particular on singular inner functions (\cite{QLS}).\end{remark}

 \begin{remark} If we do not incorporate the maximal sifting process as in UWAFD, the algorithm falls into UWA, as first developed in \cite{Na} 2000. Below we denote by $\phi_k$ Blaschke products, with finite or infinite zeros, denote by $\psi_k$ products of a singular inner function and an outer function, and denote by $f_k$  $H^2$ functions, where $f=f_1, f_k(z)=\psi_{k-1}(z)-\psi_{k-1}(0), k=2,\cdots, c_k=\psi_k(0), k=1,2,\cdots.$ It proceeds as
 \begin{eqnarray*}
 f(z)=f_1(z)&=&\phi_1(z)\psi_1(z)\\
 &=&\phi_1(z)\left(\psi_1(z)-\psi_1(0)+\psi_1(0)\right)\\
 &=&c_1\phi_1(z)+\phi_1(z)f_2(z)\\
 &=&c_1\phi_1(z)+\phi_1(z)\phi_2(z)\left(\psi_2(z)-\psi_2(0)+\psi_1(0)\right)\\
 &=&c_1\phi_1(z)+c_2\phi_1(z)\phi_2(z)+\phi_1(z)\phi_2(z)f_3(z)\\
 &=&\cdots\\
 &=&\sum_{k=1}^\infty c_k\phi_1(z)\cdots\phi_k(z).\end{eqnarray*}
The convergence in $H^2$ was first proved in Remark 4.4, \cite{Q}, and several generalized convergence results were proved in \cite{CS}. In the recent paper \cite{TQ0} we studied the computation aspect of UWA.
 \end{remark}

\subsection{Cyclic AFD for $n$-Best Rational Approximation}

In Core-AFD the parameters $a_1,\cdots ,a_k,\cdots $ are selected one by one to construct an optimal sequence of Blaschke forms to approximate the given function
\begin{eqnarray}\label{nForm} \sum_{k=1}^n \langle f,B_{\{a_1,\cdots,a_k\}}\rangle B_{\{a_1,\cdots,a_k\}}(z).\end{eqnarray}
Now we change the question to the following: Given $f\in H^2({\bf D})$ and a fixed positive integer $n,$ find $n$ parameters $a_1,...,a_n$ such that the associated $n$-Blaschke form \ref{nForm} best approximate $f,$ that is
\begin{eqnarray}\label{mim} & &\| f-\sum_{k=1}^n \langle f,B_{\{a_1,\cdots,a_k\}}\rangle B_{\{a_1,\cdots,a_k\}}(z)\|\\
&=&\min\{\|f-\sum_{k=1}^n \langle f,B_{\{b_1,\cdots,b_k\}}\rangle B_{\{b_1,\cdots,b_k\}}(z)\| \ :\ \{b_1,\cdots,b_n\}\in {\bf D}^n\}.\end{eqnarray}  This amounts an optimal but simultaneous selection of $n$ parameters that is obviously better than selections in the one by one manner.  Simultaneous selection of the parameters in an approximating $n$-Blaschke form is equivalent with the so called \emph{optimal approximation by rational functions of order not larger than} $n.$  The latter problem was phrased as \emph{$n$-best rational approximation}. It has been a long standing open problem, presented as follows.

Let $p$ and $q$ denote polynomials of one complex variable. We say that $(p,q)$ is an $n$-\emph{pair} if  $p$ and $q$ are co-prime, both of degrees less than or equal to $n,$ and $q$ does not have zero in the unit disc. Denote  by ${\cal R}_n$ the set of all $n$-pairs. If $(p,q)\in {\cal R}_n,$ then the rational function $p/q$ is said to be a rational function of degree less or equal $n.$ Let $f$ be a function in the Hardy $H^2$ space in the unit disc. To find an $n$-best rational approximation to $f$ is to find an $n$-pair $p_1,q_1)$  such that
\[  \| f-p_1/q_1 \|=\min \{ \| f-p/q \| \ : \  (p,q)\in {\cal R}_n \}.\]
Existence of such a minimum solution was proved many decades ago (\cite{Wa}), a practical algorithm to get a solution, however, has been an open problem till now. The best $n$-Blaschke form approximation is essentially equivalent with the $n$-best rational approximation. There is a separate proof for existence of the minimum in (\ref{mim})(\cite{QWe}). We wish to take the advantages of the Blaschke form to get an practical algorithm for the classical $n$-best rational approximation problem.
 By using \emph{cyclic }AFD \emph{algorithm} we can easily get a solution of the above mentioned problem if there is only one critical point for the objective function (\cite{Qcyclic}). In general, cyclic AFD offers a \emph{conditional solution} depending on the initial values to star with. Besides cyclic AFD there previously existed an algorithm, RARL2, by the French institute INRIA, that again can only get a conditional solution \cite{BCO}. Both the theory and algorithm of cyclic AFD are explicit. It directly finds out the poles of the approximating rational function. The other rational approximation models mostly use the coefficients of $p$ and $q$ as parameters in order to set up and solve the problem. Using coefficients of polynomials, which is double amount of the parameter number of the best-n Blaschke form setting, involves tedious analysis and computation. The ultimate solution of the optimization problem lays on optimal selection of an initial status to start with. Finding an optimal initial status itself is an NP hard problem.

For any given natural number $n$ the objective function for the $n$-Blaschke optimization problem is
\begin{eqnarray}\label{minimize} A(f; a_1,...,a_n)=\| f\|^2 -\sum_{k=1}^n | \langle f, B_k\rangle |^2.\end{eqnarray}\\

\noindent{\bf Definition 2}\label{Definition 1} \emph{An $n$-tuple
$(a_1,...,a_n)$ is said to be a coordinate-minimum point (CMP) of the objective function
$A(f;z_1,...,z_n)$ if for any chosen $k$ among 1,...,n, whenever we fix the rest $n-1$
variables, being $z_1=a_1,..., z_{k-1}=a_{k-1}, z_{k+1}=a_{k+1},...,z_n=a_n,$
and select the  $k$th variable $z_k$ to minimize the objective function, we have
$$a_k=\arg \min \{ A(f;a_1,...,a_{k-1}, z_k, a_{k+1},...,a_n)\ :\ z_k\in {\bf D}\}.$$ }
In the Core-AFD algorithm we progress the following procedure: For a $(k-1)$-tuple $\{a_1,...,a_{k-1}\}$ in ${\bf D}$ we produce the reduced
remainders $f_2,...,f_k,$ and for $f_k$ we apply the Maximal
Selection Principle to find an $a_k$ giving rise to max$\{
|\langle f_n, e_{a}\rangle | : a\in  {\bf D}\}.$ The proposed cyclic AFD
algorithm repeats such procedure always for $k=n:$ For any permutation $P$ of
$1,...,n,$ for the first $(n-1)$ parameters in the order
$a_{P(1)},...,a_{P(n-1)}$ we produce the corresponding reduced remainders $f_2,...,f_n,$ and then use the
Maximal Selection Principle to select a new and optimal $a_{P(n)}.$

The proposed cyclic AFD algorithm is validated in the following theorem.

\begin{theorem} Suppose that $f$ is not an $m$-Blaschke
form for any $m<n.$ Let $s_0=\{b_1^{(0)},...,b_n^{(0)}\}$ be any
$n$-tuple of parameters inside ${\bf D}.$ Fix some $n-1$
parameters of $s_0$ and make an optimal selection of the single
remaining parameter according to the Maximal Selection Principle based on the objective function
(\ref{minimize}).
Denote the obtained new $n$-tuple of parameters by $s_1.$ We
repeat this process and make cyclic optimal selections over the
$n$ parameters. We thus obtain a sequence of $n$-tuples $s_0, s_1,
..., s_l,...,$ with decreasing objective function values $d_l$
that tend to a limit $d\ge 0,$ where, in the notation and
formulation of (\ref{minimize}),
\begin{eqnarray}\label{DL}
d_l=A(f; b_1^{(l)},...,b_n^{(l)})=\| f\|^2- \sum_{k=1}^n
(1-|b_k^{(l)}|^2) |f_k^{(l)} (b_k^{(l)})|^2.
\end{eqnarray}  Then, (i) If $\overline{s},$ as an $n$-tuple, is
a limit of a subsequence of $\{s_l\}_{l=0}^\infty,$ then
$\overline{s}$ is in ${\bf D};$ (ii) $\overline{s}$ is a CMP of
$A(f;\cdots );$ (iii) If the correspondence between a CMP and the
corresponding value of $A(f;\cdots )$ is one to one, then the
sequence $\{s_l\}_{l=0}^\infty$ itself converges to the CMP, being
dependent of the initial $n$-tuple $s_0;$ (iv) If $A(f;...)$ has
only one CMP, then $\{s_l\}_{l=0}^\infty$
converges to a limit $\overline{s}$ in ${\bf D}$ at which $A(f;\cdots )$
attains its global minimum value.\end{theorem}

We refer the reader to \cite{Qcyclic} for further details and examples of cyclic AFD. In a recent paper the algorithm is further developed \cite{QWj}.

\subsection{Pre-Orthogonal adaptive Fourier Decomposition (POAFD) for Reproducing Kernel Hilbert Spaces}

The theory and algorithm that will be developed in this section, as a matter of fact, can be extended to more general contexts. To explain just the idea we restrict ourselves to the simplest cases, including the weighted Bergman spaces and weighted Hardy spaces, etc. In the said simple setting the Hilbert space ${\cal H}$ consists of functions defined in an open connected region ${\cal E}$ (can be unbounded) in the complex plane, and the reproducing kernel $k_a$ is a real-analytic function of the variable $\overline{a}$ in ${\cal E}$ satisfying the relation
\begin{eqnarray}\label{partial} f^{(l)}(a)=\langle f,\left(\frac{\partial}{\partial \overline{a}}\right)^{l}k_a\rangle, \quad l=1,2,\cdots .\end{eqnarray}
Let $\{a_1,\cdots,a_n,\cdots\}$ be a finite or infinite sequence. For a fixed $n$ we define the multiple of $a_n,$ denoted by $l(a_n),$ to be the repeating times of $a_n$ in the $n$-tuple $\{a_1,\cdots,a_n\}.$ With this definition, for instance, the multiple of $a_1$ is just $1,$ and the multiple of $a_2$ will depend on whether $a_2=a_1.$ If yes, then $l(a_2)=2,$ and, if not, $l(a_2)=1,$ and so on. Note that it is a little abuse of notation for it is not dependent on the value of $a_n$ but on the repeating times of $a_n$ in the corresponding $n$-tuple. We accordingly define \begin{eqnarray}\label{tilde}\tilde{k}_{a_n}\triangleq
\left[\left(\frac{\partial}{\partial\overline{a}}\right)^{l(a_n)-1}
k_{a}\right]_{a=a_n}\triangleq\left(\frac{\partial}{\partial\overline{a}}\right)^{l(a_n)-1}
k_{a_n}.\end{eqnarray} We further assume the following boundary vanishing condition, implying the Maximal Selection Principle in every individual context, as follows: Let $a_1,\cdots,a_{n-1}$ be previously given, and $\{B_1,\cdots,B_{n-1}\}$ be the Gram-Schmidt orthonormalization of $\{\tilde{k}_{a_1},\cdots,\tilde{k}_{a_{n-1}}\},$ then for every $f\in {\cal H},$ the pre-orthogonal system has the property
\begin{eqnarray}\label{boundary vanishing} \lim_{a\to \partial{\cal H}}\langle f, B^a_n\rangle =0,\end{eqnarray}
where $\{B_1,\cdots,B_{n-1},B^a_n\}$ is the Gram-Schmidt orthonormalization of $\{\tilde{k}_{a_1},\cdots,\tilde{k}_{a_{n-1}}, k_a\},$ with $a\ne a_k, k=1,\cdots, n-1.$ We note (1) if $a\to \partial {\cal H}$ then $a$ is different from any $a_k, k=1,\cdots,n-1,$ when $a$ gets close to the boundary; and (2) in any case the limit $a\to \partial {\cal H}$ is in the sense of the topology of the one-point-compactification of the complex plane while the $\lq\lq$one point" takes to be $\infty.$ With this boundary vanishing assumption we conclude the Maximal Selection Principle of POAFD: Under the assumption (\ref{boundary vanishing}), through a compact argument using the Bolzano-Weierstrass theorem, there exists a sequence $\{b_j\}_{j=1}^\infty$ such that none of the $b_j$'s take any values $a_1,\cdots,a_{n-1},$ and $\lim_{j\to \infty}b_j=a_n\in {\cal E},$ and
\begin{eqnarray}\label{maximal}
\lim_{j\to\infty}|\langle f,B^{b_j}_n\rangle|=\max\{|\langle f,B^a_n\rangle|\ :\ a\in {\cal H}\}.\end{eqnarray}

Under those conditions we can prove the following lemma.

\begin{lemma}
\[\lim_{l\to\infty}B^{b_j}_n=B^{a_n}_n,\]
where $\{B_1,\cdots,B_{n-1},B^{a_n}_n\}$ is the Gram-Schmidt orthonormalization of $\{\tilde{k}_{a_1},\cdots,\tilde{k}_{a_{n-1}},\tilde{k}_{a_{n}}\}.$
\end{lemma}

\noindent{\bf Proof}. If $a_n$ does not coincide with any $a_k, k=1,\cdots, a_{n-1},$ then $\lim_{j\to\infty}B^{b_j}_n=B^{a_n}_n,$ where $\{B_1,\cdots,B_{n-1},B^{a_n}_n\}$ is the Gram-Schmidt orthonormalization of $\{\tilde{k}_{a_1},\cdots,\tilde{k}_{a_{n-1}}, k_{a_{n}}\}=\{\tilde{k}_{a_1},\cdots,\tilde{k}_{a_{n-1}},\tilde{k}_{a_{n}}\}.$
Now consider the case that $a_n$ coincides with some of the earlier $a_1,\cdots,a_{n-1},$ or in other words, $l(a_n)>1.$
That means that, in the notation (\ref{tilde}), the $(l-1)$ functions $k_{a_{n}}, \frac{\partial}{\partial \overline{a}}k_{a_{n}},\cdots, $ $\left(\frac{\partial}{\partial \overline{a}}\right)^{(l-2)}k_{a_{n}}$ have already appeared in the sequence $\{\tilde{k}_{a_1}, \cdots,\tilde{k}_{a_{n-1}}\}.$ As a consequence, the function
\[ T_{l-2}(b_j, a_{n})=k_{a_n}+\frac{\frac{\partial}{\partial \overline{a}}k_{a_{n}}}{1!}(\overline{b_j}-\overline{a}_{n})+\cdots
+\frac{\left(\frac{\partial}{\partial \overline{a}}\right)^{(l-2)}k_{a_{n}}}{(l-2)!}(\overline{b_j}-\overline{a}_{n})^{l-2},\]
as the order-$(l-2)$ Taylor expansion of the function $k_a(z)$ in $\overline{b}_j$ about $\overline{a}_{n},$ is already in the linear span of $B_1,\cdots, B_{n-1}.$ This last atatement amounts to the relation
 \begin{eqnarray}\label{insert1} T_{l-2}(b_j,a_{n})-\sum_{k=1}^n\langle T_{l-2}(b_j,a_{n}),B_k\rangle B_k=0.\end{eqnarray}
 Since $b_j$'s are different from all $a_k, k=1,\cdots, n,$ we have,
\begin{eqnarray}\label{into1} B_{n}^{b_j}(z)=\frac{k_{b_j}(z)-\sum_{k=1}^{n-1} \langle k_{b_j}, B_k\rangle B_k (z)}{\|k_{b_j}-\sum_{k=1}^{n-1} \langle k_{b_j}, B_k\rangle B_k \|},\end{eqnarray}

Inserting (\ref{insert1}) into (\ref{into1}), and dividing by $(\overline{b_j}-\overline{a}_{n})^{l-1}$ and
$|\overline{b_j}-\overline{a}_{n}|^{l-1}$ to the numerator and the denominator parts, respectively, we have
\begin{eqnarray}\label{inserted} B_{n}^{b_j}(z)=e^{-(l-1)\theta}\frac{\frac{k_{b_j}(z)-T_{l-2}(k_{b_j}, a_{n})(z)}{(\overline{b_j}-\overline{a}_{n})^{l-1}}-\sum_{k=1}^{n-1} \langle \frac{k_{b_j}-T_{l-2}(b_j, a_{n})(z)}{\overline{w}-\overline{a}_{n})^{l-1}}, B_k\rangle B_k (z)}{\|\frac{k_{b_j}(z)-T_{l-2}(b_j, a_{n})(z)}{(\overline{b_j}-\overline{a}_{n})^{l-1}}-\sum_{k=1}^{n-1} \langle \frac{k_{b_j}(z)-T_{l-2}(b_j, a_{n})(z)}{(\overline{b_j}-\overline{a}_{n})^{l-1}}, B_k\rangle B_k \|},\end{eqnarray}
where $e^{i\theta}$ is the tangential direction of the limiting $b_j\to a_n.$ We can, in fact, take any direction, including $\theta=0.$
Letting $b_j\to a_n$ with $\theta=0,$ and using the Lagrange type remainder of the Taylor expansion, we obtain
\[ \lim_{j\to \infty}B^{b_j}_n(z)=\frac{\tilde{k}_{a_n}(z)-\sum_{k=1}^{n-1} \langle \tilde{k}_{a_n}, B_k\rangle B_k (z)}{\|\tilde{k}_{a_n}-\sum_{k=1}^{n-1} \langle \tilde{k}_{a_n}, B_k\rangle B_k \|}.\]
Therefore, $\{B_1,\cdots,B_{n-1},B^{a_n}_n\}$ is the Gram-Schmidt orthonormalization of $\{\tilde{k}_{a_1},\cdots,\tilde{k}_{a_{n-1}},\tilde{k}_{a_{n}}\},$ as desired.

\begin{remark}
The essence of the proof is contained in \cite{Q2D} and \cite{Qbook}. The proof for the one complex variable case as presented here first appears in \cite{QD}.
\end{remark}

We have the pre-orthogonal adaptive Fourier decomposition (POAFD) convergence theorem as follows.

\begin{theorem}
Selecting $\{a_1,\cdots,a_n,\cdots\}$ according to the Maximal Selection Principle set by (\ref{maximal}), we have
\[ f=\sum_{k=1}^\infty \langle f,B_n\rangle B_n,\]
where for any positive integer $n,$ $\{B_1,\cdots,B_{n-1},B_n\}$ is the Gram-Schmidt orthonormalization of $\{\tilde{k}_{a_1},\cdots,\tilde{k}_{a_{n-1}},\tilde{k}_{a_{n}}\}.$
\end{theorem}
One can adopt the same proof for the AFD convergence (Theorem \ref{AFD}) in which only the standard remainders $g_k$'s are concerned. As a matter of fact, the sifting process and the role of the induced remainders are taken place by the pre-orthogonal process.
\begin{remark}For Repeating selection of parameters and POAFD we refer the reader to the references \cite{QWe,QSW,Q2D,Qbook}. was previously called POGA or PreOGA, etc. The previous names do not reflect the crucial role of the complex holomorphic function methods.\end{remark}
\section{Related Studies and Applications}

\subsection{Aspects in Relation to Beurling-Lax Shift-invariant Subspaces}

The AFD type expansions is, in a great extent, related to the Beurling-Lax shift-invariant subspaces of the Hardy $H^2$ spaces. In the unit disc case,
\begin{eqnarray}\label{direct sum} H^2({\bf D})=\overline{\rm span}\{B_k\}_{k=1}^\infty\oplus \phi H^2({\bf D}),\end{eqnarray}
where $\{B_k\}_{k=1}^\infty$ is the TM system generated by a sequence $\{a_1,\cdots,a_n,\cdots\},$ where multiples are counted, and $\phi$ is the Blaschke product with the zeros $\{a_1,\cdots,a_n,\cdots\}$ including the multiples. Note that when $\phi$ can be defined with the $a_k$'s as its zeros, there holds the condition
\[ \sum_{k=1}^\infty (1-|a_k|)<\infty,\] and  the associated TM system is not a basis. Although this has been well know its relations with adaptive expansions, as far as being aware by the author, are for the first time explored. Lately TM systems being Schauder system were proved (\cite{QCT}). The space decomposition relation (\ref{direct sum}) was extended to $H^p$ spaces, where $p\ne 2.$ Relations between backward shift invariant subspaces and bandlimited functions and Bedrosian identity (\cite{QT2}) were studied. There are open problems such as whether for $p\ne 2$ there exist adaptive and fast converging expansions by using TM systems, and for $p=2$ how far one can extend (\ref{direct sum}) to higher dimensions. The study has a great room to be further developed.

\subsection{Extra-strong Uncertainty Principle}

The phase and frequency studies in the mono-component function theory lay necessary foundations for digital signal processing. In the related studies what is called extra-strong uncertainty principle
\begin{eqnarray}\label{extra} \sigma_t^2\sigma_\omega^2\ge \frac{1}{4}+\left(\int_{-\infty}^\infty
|t-\langle t\rangle||\phi (t)-\langle \omega \rangle||f(t)|^2dt\right)^2,\end{eqnarray}
was established, where $f$ is a real-valued signal, $\sigma_t^2$ and $\sigma_\omega^2$ are the standard deviations with respect to, respectively, the time and the Fourier frequency, and  $\langle t\rangle $ and $\langle \omega \rangle$ are the corresponding means (\cite{DDQ1}). An uncertainty principle of the same type was given by L. Cohen
\[ \sigma_t^2\sigma_\omega^2\ge \frac{1}{4}+|\int_{-\infty}^\infty
(t-\langle t\rangle)(\phi (t)-\langle \omega \rangle)|f(t)|^2dt|^2\]
which is obliviously weaker. We further extend the above result to multi-dimensional contexts \cite{DDQ1,DDQ2,DQY,DQY1,DQC}.

\subsection{The Dirac-Type Time-Frequency Distributions Based on Mono-component Decompositions}
The Dirac type time-frequency distribution (D-TFD) of the form
\begin{eqnarray}\label{DTFD} P(t,\omega)=\rho^2(t)\delta (\omega-\theta'(t))\end{eqnarray}
 has been ultimate destiny of signal analysts. There have existed several time-frequency distributions, including windowed Fourier transform and Wagner-Ville transform, etc., of which none are satisfied. The existing time-frequency distributions do not give explicit and clear frequency components, and often depend on parameter selections. Positive-frequency decompositions of signals offered by the AFD decompositions naturally give rise to Dirac-type time-frequency distributions. For a single mono-component $m_1(t)=\rho_1(t)\cos\theta_1(t)$ the corresponding D-TFD according to (\ref{DTFD}) is the graph $(t,\theta_1'(t))$ of the function $\omega =\theta_1'(t)$ in the $\omega$-$t$ plane, while the weight $\rho_1^2(t)$ may be represented by colors continuously changing along with changing of the values $\rho_1^2(t).$ If a signal $f$ is expanded into a series consisting of its $\lq\lq$intrinsic composing" mono-components, then its D-TFD is the bunch of color-weighted graphs of which each is made from a composing mono-component (\cite{ZQDM,DQ4}). This definition has been interested and being paid attention by prominent signal analysts including Leon Cohen and Lorenzo Galleani, etc., and has been used in practice (see below the application section).
\subsection{Higher Dimensional AFDs}
To develop an AFD like approximation theory in higher dimensions a Cauchy type structure is necessary, that is mainly for use of the reproducing kernel property in deducing Maximal Selection Principle in the underlying space. By using the Cauchy structure in Clifford algebra or in several complex variables we achieved the AFD type theories and algorithms for functions of several real variables on the plane (Clifford Hardy spaces and Hardy spaces on tubes), and on the real spheres, and of several complex variables on the $n$-torus, and on the $n$-complex spheres (\cite{QSW,QWY,WQ1,QWY,SQSW}). With D. Alpay, F. Colombo, I. Sabadini we achieved analogous theory involving matrix valued Blaschke products (\cite{ACQS1,ACQS2}. This study has impact in general to rational approximation in a number of spaces (\cite{BMQ}).

\subsection{Fourier Spectrum Characterization of Hardy Spaces: Analytic Signals Revised}

The Paley-Wiener Theorem for the classical Hardy $H^2$ space over the upper-half complex space addresses the fact that if $f\in L^2({\bf R}),$ then further $f\in H^2({\bf C}^+)$ if and only if ${\rm supp}\hat{f}\subset [0,\infty).$  This result is systematically extended to  $H^p({\bf C}^+)$ for all $p\in [1,\infty],$ where the Fourier transform can be defined at least in the distribution sense (\cite{Q17,QXYYY}). The results are summarized as: Letting $f\in L^p({\bf R}),$ then $f$ is further the non-tangential boundary limit of some function in the complex Hardy space $H^p({\bf C}^+)$ if and only if $\hat{f}=\chi_+\hat{f},$ where is the indicator (characteristic) function of the right-half-real line, and the Fourier transform may take the distribution sense. The generalization to the Hardy spaces on tubes (extending the $p=2$ case in \cite{SW} to $1\leq p\leq \infty$) was published in \cite{LDQ1}. The results of the same type but with the Clifford algebra setting is now on line \cite{DMQ1} proving that a Clifford-valued function $f\in L^p({\bf R}^n)$ is the non-tangential boundary limit of some Clifford-valued Hardy space function in the upper-half space if and only if $\hat{f}=\chi_+\hat{f},$ where $\chi_+(\underline{\xi})=\frac{1}{2}\left( 1+i\frac{\underline{\xi}}{|\underline{\xi}|}\right), \underline{\xi}=\xi_1e_1+\cdots +\xi_ne_n$ (the Hardy space projection function). This extends some partial cases proved in \cite{SW} for conjugate harmonic systems. In various contexts Fourier spectrum characterizations give rise to Hardy spaces decompositions for $L^p, \ 1<p\leq\infty$ that further induce Hardy space decompositions of $L^p,\ 0<p<1$ (\cite{DeQ,LDQ2}). Hardy space decomposition is the strategy that we have been using to study Lebesgue spaces of various integrability. The strategy is extensively implemented along with the mono-component and AFD approximation theories. In particular, for any signal $f$, by multiplying $\hat{f}$ with $\chi_+$ and then taking the inverse Fourier transform, we obtain the associated analytic signal. This is philosophically valid in any context. We finally note that the Hardy space decomposition issue has been extended to the $L^p$-vector fields and one obtains the Hardy-Hodge decomposition \cite{BDQ}.

\subsection{Hilbert Transforms as Singular Integral Operators: Analytic Signals Revised}

As in the 1-D case (\cite{Be}), on higher dimensional manifolds one defines the non-scalar part of the non-tangential boundary limit of a hyper-complex holomorphic function to be the Hilbert transform of the scalar part of it (\cite{QY}). Hilbert transform therefore is a particular singular integral. It is, in particular, not the singular Cauchy transform. One must study singular integrals to understand Hilbert transform. On one dimensional manifolds, including Lipechitz perturbations of the real line and the circle, certain singular integrals of holomorphic kernels form an operator algebra as studied in a series of work of A. McIntosh, C. Li, S. Semmes, T. Qian, R.-L. Long and S.-L. Wang (\cite{McQ1,McQ2,Q11,GQW}). The theory on the plane was earlier established in the work or under the influence of A. McIntosh (\cite{LMcQ, LMcS, GLQ}). Through first generalizing the results of Fueter and Sce to arbitrary Euclidean spaces (as technical necessity) the author established the theory of the operator algebra of singular integrals of Clifford monogenic kernels on Lipschitz perturbations of the unit sphere for any dimension (\cite{Q15}). Based on the established singular integral theory Hilbert transformations of the plane and of the sphere became well understood. Analytic signals on the sphere, for instance, are constructed as follows. Let $f$ be a real-valued signal of finite energy on a manifold ${\cal S}.$  Denote by $H_{\cal S}$ the Hilbert transform of $f$ on the manifold. Then the analytic signal on ${\cal S}$ is defined to be $f^+=f+H_{\cal S}f,$ where $H_{\cal S}f$ is the non-scalar part (on sphere it is a 2-form valued function). $f^+$ has Clifford monogenic extension to one of the two regions divided by ${\cal S}.$ One can derive, if $\zeta$ is on the plane or on the sphere,
\begin{eqnarray}\label{manifold}
f^+(\zeta)&=& \rho_f(\zeta)\left(\frac{f(\zeta)}{\rho_f(\zeta)}+\frac{H_{\cal S}(f)(\zeta)}{|H_{\cal S}(f)(\zeta)|}\frac{|H_{\cal S}(f)(\zeta)|}{\rho_f (\zeta)} \right)\\
&=&\rho_f (\zeta)\left(\cos\theta(\zeta)+\frac{H_{\cal S}(f)(\zeta)}{|H_{\cal S}(f)(\zeta)|}\sin\theta(\zeta) \right)\\
&=& \rho_f (\zeta) e^{\frac{H_{\cal S}(f)(\zeta)}{|H_{\cal S}(f)(\zeta)|}\theta (\zeta)},\end{eqnarray}
where $\rho_f (\zeta)=\sqrt{|f(\zeta)|^2+|H_{\cal S}f(\zeta)|^2},$ and $\left(\frac{H_{\cal S}(f)(\zeta)}{|H_{\cal S}(f)(\zeta)|}\right)^2=-1,$
the latter being a varying imaginary element just like the complex imaginary element with the property $i^2=-1.$ The instantaneous frequency is defined, as in the classical case through the monogenic continuation, but formally read
\[ \theta'(\zeta)={\rm Re}\{\left[ (\Gamma_\zeta -I)f^+(\zeta)\right]\left[(f^+(\zeta))^{-1}\right]\},\] the latter can be expressed in terms of the angle $\theta(\zeta)$ {\cite{QY, YQS2}, where $\Gamma_\zeta$ is the surface Dirac operator on the manifold. In such format the basic idea of IF and the related approximation in higher dimensions make sense. The related studies published, in respectively 2015 and 2017 with the Chinese Science Press two monographs books \cite{Qbook} and \cite{Qbook2}. We finally note that Hilbert transformation may be characterize by commutativity with the affine groups in the underlying symmetric manifold which shows that the three objects the Hilbert transformation, the Dirac differential operator and the group representation theory have intimate relations (\cite{DLQ1,DLQ2}).
\subsection{Applications}

AFD has demonstrative applications in system identification and signal analysis. Applications in system identification include \cite{MQ1,MQW,MQM,CMZM}). A number of signal analysts promoted the AFD methods. Below we summarize part of the applications found in the literature.

It is commented in \cite{WW1} that as a new method AFD was proposed in the recent years that could be used to decompose and reconstruct signals. It contains the classical Fourier method as a particular case. Experiments show that the 1D AFDs achieve excellent signal decomposition and reconstruction results. The article [W1] compares 2D AFD with the traditional frequency digital watermark methods, including discrete cosine transform DCT, discrete wavelet transform DWT, discrete Fourier transform DFT, etc., and concludes that 2D AFD has better transparency and robustness under attacking. The articW2] revises the 2D AFD algorithm and, as a result, increases its speed, and uses it in denoising.

In \cite{LJC} Y. Liang et al. at Beijing Jiaotong University propose a new fault diagnosis method of rolling bearing based on AFD. They show that AFD can avoid using band-pass filters, the latter often suffering from the difficulty of algorithm parameter selection, they show that AFD adaptively, efficiently, and accurately diagnose all kinds of rolling bearing problems.

In \cite{WCW} the authors study interference and separation between the lung sound (LS) and the heart sound (HS) signals. Due to the overlap in their frequency spectra, it is difficult to separate them. The article proposes a novel separation method based on AFD. This AFD-based separation method is validated on real HS signals from the University of Michigan Heart Sound and Murmur Library, as well as real LS signals from the 3M repository. Simulation results indicate that the proposed method is more effective than the extraction methods based on the recursive least square (RLS), than the standard empirical mode decomposition (EMD) and its various extensions, including the ensemble EMD (EEMD), the multivariate EMD (M-EMD) and the noise assisted M-EMD (NAM-EMD).

Over the years people have made unremitting studies in predicting the stock price movements. In \cite{Z1} a novel automatic stock movement forecasting system is proposed, which is based on the newly developed signal decomposition approach - adaptive Fourier decomposition (AFD). AFD can effectively extract the signal primary trend, which is specifically suitable in the Dow Theory based automatic technique analysis. Effectiveness of the proposed approach is assessed through the comparison with the direct BP approach and manual observation. The result is proved to be promising.

In \cite{Z2} an AFD based time-frequency speech analysis approach is proposed. Given the fact that the fundamental frequency of speech signals often undergo fluctuation, the classical short-time Fourier transform (STFT) based spectrogram analysis suffers from the difficulty of window size selection. AFD is a newly developed signal decomposition theory. The outstanding characteristic of AFD is to provide instantaneous frequency for each decomposed component, so the time-frequency analysis becomes accessible. Experiments are conducted based on the sample sentence in TIMIT Acoustic-Phonetic Continuous Speech Corpus. The results show that the AFD based time-frequency distribution outperforms the STFT.

AFD has already been employed to the productions of CASA Environmental Technology Co., Ltd, including the second generation of BEWs (Biological Early Warning System) and ETBEs (Ecological Toxicity Biological Exposed System). The two systems takes advantages of AFD and Unwinding to analyze the biological behavioral signal. Compared with the traditional Fourier, wavelet and EMD algorithms, the AFD approach can efficiently solve early warning judgment for low concentration pollutants and disturbance of fish biological clock and other problems.

In control theory the the authors of the article \cite{LFZ} introduce an AFD algorithm to eliminate the channel noise superimposed on the output signal in the wireless transmission process. In the frequency domain, based on AFD, an ILC method for discrete linear system with wireless transmission is proposed. Simulation results show that the AFD algorithm is able to achieve signal denoising well in the case of small decomposition threshold compared with Fourier decomposition. Thus the goal that the output signal of ILC system can track the desired signal is better achieved.

Indian researchers in their article \cite{GMA} assert that to analyze biomedical signals in relation to e-health devices the frequency domain method outperforms the time domain method, and among numerate frequency domain methods (Hermit, Fourier, Karhunen-Loeve, Wavelet) AFD appears to have features of a greater variety, and more stable for the data compression. Based on compression using AFD they started to manufacture economic, accurate and stable domestic e-health devices.

Apart from China and Asia, AFD has also achieved international influence. Interests, studies and applications of AFD are found in relevant literature, by Ph.D. thesis of F. D. Fulle at Michigan University on oxygenic photosynthesis; by A. Kirkbas et al. on optimal basis pursuit based on jaya optimization for adaptive Fourier decomposition (\cite{KKB}); by V. Vatchev, on a class of intrinsic trigonometric mode polynomials (\cite{VVV}); by J. Mashreghi et al. on Blaschke Products and Applications (\cite{MF}); by R.S. Krausshar et al. on Clifford and harmonic analysis on cylinders and tori (\cite{KR}); by F. Colombo et al. on the Fueter mapping theorem in integral form and the ${\cal F}$-functional calculus (\cite{CSS}); by M.I. Falc$\tilde{a}$o et al. on remarks on the generation of monogenic functions; by F. Colombo et al. on the Fueter primitive of bi-axially monogenic functions (\cite{CSS0}); by L. Salomon on analysis of the anisotropy in image textures (\cite{SA}); by F. Sakaguchi on the related integral-type method in higher order differential equations (\cite{SH1,SH2,SH3,SH4}); by P. Le\'on on instantaneous frequency estimation and representation of the audio signal through complex wavelet additive synthesis (\cite{LBB}); by F.E. Mozes on computing the instantaneous frequency for an ECG signal (\cite{MS}); by N.R. Gomes, as Doctoral dissertation, on compressive sensing in Clifford analysis; by T. Eisner et al. on discrete orthogonality of the Malmquist Takenaka system on the upper half plane and rational approximation (\cite{EP}); and by A. Perotti on his article in directional quaternionic Hilbert operators (\cite{Pe}).


\begin{thebibliography}{AFD}

\bibitem{ACQS1} D. Alpay, F. Colombo, T. Qian, I. Sabadini, {\it Adaptive orthonormal systems for matrix-valued functions,} Proceedings of the American Mathematical Society, 2017, 145(5)�G2089�V2106.
 \bibitem{ACQS2} D. Alpay,  F. Colombo, T. Qian, and I. Sabadini, {\it Adaptive Decomposition: The Case of the Drury-Arveson Space,} Journal of Fourier Analysis and Applications, 2017, 23(6), 1426-1444.

        \bibitem{AKQ}  A. Axelsson, K. I. Kou, T. Qian, {\it Hilbert transforms and the Cauchy integral in Euclidean space,} Studia Mathematica, 2009, 193(2): 161-187.


\bibitem{BCO}  L. Baratchart, M. Cardelli, M. Olivi, Identification and rational $L^2$ approximation,
a gradient algorithm, Automatica, 27(1991), pp. 413-418.

\bibitem{BDQ} L. Baratchart, P. Dang, T. Qian, {\it Hardy-Hodge Decomposition of Vector Fields in Rn,} Transactions of the American Mathematical Society, 2017, 1-19.

\bibitem{BMQ} L. Baratchart, W.X. Mai, T. Qian, {\it Greedy Algorithms and Rational Approximation in One and Several Variables,} In: Bernstein S., Kaehler U., Sabadini I., Sommen F. (eds) Modern Trends in Hypercomplex Analysis. Trends in Mathematics, pp: 19-33, 2016.

\bibitem{Be} S. Bell, {\it The Cauchy Transform, Potential theory and Conformal Mappings,} CRC Press, Boca, Raton (1992).

\bibitem{Boa} B. Boashash, {\it Estimating and interpreting the instantaneous frequency of
a signal-Part 1: Fundamentals,} Proceedings of The IEEE, vol.80, no.4,
pp.520-538, 1992.

\bibitem{CD} M. T. Cheng and D. G. Deng, {\it Lecture notes on harmonic analysis,} Beijing University, 1979.

\bibitem{CP} R. Coifman and J. Peyri\'ere, {\it Phase unwinding, or invariant subspace decompositions of Hardy spaces,} arXiv.org > math > arXiv:1707.04844v1.

\bibitem{CS} R. Coifman R, S. Steinerberger, {\it Nonlinear phase unwinding of functions,} J Fourier Anal Appl, 2017, 23: 778��?09.

    \bibitem{CSW} R. Coifman, S. Steinerberger, H-t. Wu, Carrier frequencies, holomorphy and unwinding, SIAM J. Math. Anal., accepted.


\bibitem{ChengQS} Q.S. Cheng, {\it Digital Signal Processing,} Peking University Press,
2003, in Chinese.

\bibitem{ChenQTan} Q-H. Chen, T. Qian and L-H. Tan, {\it Constructive Proof of Beurling-Lax Theorem,} accepted to appear in Chin. Ann. of Math.

\bibitem{Co} L. Cohen, {\it Time-Frequency Analysis: Theory and Applications,} Prentice Hall, 1995.

\bibitem{CMZM} Q.-H. Chen, W.-X.Mai, L.-M. Zhang, W. Mi, {\it  System identification by discrete rational atoms,} Automatica, 2015, 56. pp. 53-59.

\bibitem{CSS0} F. Colombo, I. Sabadini, F. Sommen F, {\it The Fueter primitive of bi-axially monogenic functions,} Communications on Pure and Applied Analysis, 2014, 13(2).

\bibitem{CSS} F. Colombo, I. Sabadini, F. Sommen, {\it  The Fueter mapping theorem in integral form and the ${\cal F}$-functional calculus,} Mathematical Methods in the Applied Sciences, 2010, 33(17): 2050-2066.

    \bibitem{DDQ1} P. Dang, G. T. Deng, T. Qian, {\it A Sharper Uncertainty principle,} Journal of Functional Analysis, 2013, 265(10): 2239-2266.

\bibitem{DDQ2} P. Dang, G. T. Deng, T. Qian, {\it A Tighter Uncertainty Principle For Linear Canonical Transform in Terms of Phase Derivative,} IEEE Transactions on Signal Processing, 2013, 61(21): 5153 - 5164.

    \bibitem{DLQ1} P. Dang, H. Liu, T. Qian, {\it Hilbert Transformation and Representation of $ax+b$ Group,} accepted by Canadian Mathematical Bulletin, 61(1):1-15, October 2017. DOI:10.4053/CMB-2017-063-0.

    \bibitem{DLQ2} P. Dang, H. Liu T. Qian, {\it Hilbert Transformation and $r{\rm Spin}(n)+{\bf R}^n$ Group,} arXiv:1711.04519v1[math.CV]

\bibitem{DMQ1} P. Dang, W.X. Mai, T. Qian, {\it Fourier Spectrum Characterizations of Clifford $H^p$ Spaces on $R_+^{n+1}$ for $1\leq p\leq \infty,$} arXiv:1711.02610[math.CV].




\bibitem{DQ2} P. Dang, T. Qian,{\it Analytic Phase Derivatives, All-Pass Filters and Signals of Minimum Phase,} IEEE Transactions on Signal Processing, 2011, 59(10): 4708 �V 4718.


\bibitem{DQ4} P. Dang, T. Qian, {\it Transient Time-Frequency Distribution based on Mono-component Decompositions,} International Journal of Wavelets, Multiresolution and Information Processing, 2013, 11(3), 1350022, 24 pp.


 \bibitem{DQC} P. Dang, T. Qian,  Q. H. Chen, {\it Uncertainty Principle and Phase�VAmplitude Analysis of Signals on the Unit Sphere,} Advances in Applied Clifford Algebras, 2017, 27(4), 2985-3013

    \bibitem{DQY} P. Dang, T. Qian and Z. You, {\it Hardy-Sobolev spaces decomposition and applications in signal analysis,} J. Fourier Anal. Appl. 17 (2011), no. 1, 36-64.

\bibitem{DQY1} P. Dang, T. Qian, Y. Yang, {\it Extra-strong uncertainty principles in relation to phase derivative for signals in Euclidean spaces,} Journal of Mathematical Analysis and Applications, 2016, 437(2)�G912-940.

     \bibitem{DQY} P. Dang, T. Qian, Z. You, {\it Hardy-Sobolev spaces decomposition and applications in signal analysis,} Journal of Fourier Analysis and Applications,  2011, 17(1): 36�V64.

         \bibitem{DeQ} G.T Deng, T. Qian, {\it Rational approximation of Functions in Hardy Spaces,} Complex Analysis and Operator Theory, 2016, 10(5), pp. 903-920.

    \bibitem{EP} T. Eisner, M. Pap, {\it Discrete orthogonality of the Malmquist Takenaka system of the upper half plane and rational
     interpolation,} Journal of Fourier Analysis and Applications, 2014, 20(1): 1-16.


 \bibitem{FCM} M.I. Falc$\tilde{a}$o, J.F. Cruz, H.R. Malonek, {\it Remarks on the generation of monogenic functions,} International Conference on the Applications of Computer Science and Mathematics in Architecture and Civil Engineering, 17, Weimar, 2006.

\bibitem{FM} P. Fulcheri, M. Olivi, Matrix rational $H^2$ approximation: a gradient algorithm
based on schur analysis, SIAM I. Control Optim., Vol 36, No. 6, pp. 2103-2127, November 1998.

\bibitem{Gab} D. Gabor, {\it Theory of communication,} J. IEE. 93(III), pp 429-457, 1946.

\bibitem{GLQ}  G. I. Gaudry, R. Long, T. Qian, {\it A Martingale proof of L2-boundednessof Clifford-Valued Singular Integrals,} Annali di Mathematica Pura Ed Applicata, 1993, 165: 369-394.

\bibitem{GQW} G. Gaudry, T. Qian, S. L. Wang, {\it Boundedness of singular integrals with holomorphic kernels on star-shaped closed Lipschitz curves,}  Colloquium Mathematicum, 1996, LXX: 133-150.

\bibitem{Ga} J.B. Garnett, Bounded Analyic Functions, Academic
Press, 1981.

 \bibitem{Go} N.R. Gomes, {\it Compressive sensing in Clifford analysis,} Doctoral Dissertation, Universidade de Aveiro (Portugal), 2015.

\bibitem{GS} S. Gong, {\it Private comminication}, 2002.

\bibitem{Gorusin}, G.M. Gorusin, {\it Geometrical Theory of Functions of One Complex Variable,} translated by Jian-Gong Chen, 1956.

\bibitem{GMA} P. Ganta, G. Manu, and S. Anil Sooram, {\it New Perspective for Health Monitoring System,} International Journal of Ethics in Engineering and Management Education, ISSN:2348-4748, Volume 3, Issue 10, October 2016.

    \bibitem{Hummel} J.A. Hummel, {\it Multivalent starlike function,} J. d' analyse Math. 18, 133��?60 (1967)


    \bibitem{Huang} N. E. Huang et al, {\it The empirical mode decomposition and the Hilbert spectrum for
nonlinear and non-stationary time series analysis,} Proc. R. Soc. Lndon, A454(1998),
903-995.

\bibitem{KKB} A. Kirkbas, A. Kizilkaya, E. Bogar, {\it Optimal basis pursuit based on jaya optimization for adaptive fourier decomposition,} Telecommunications and Signal Processing, 2017 40th International Conference on IEEE: 538-543.

\bibitem{KR} R.S. Krausshar, J. Ryan, {\it Clifford and harmonic analysis on cylinders and tori,} Revista Matematica Iberoamericana, 2005, 21(1): 87-110.

\bibitem{LBB} P. de Le\'on, J.R. Beltr\'an, Beltr\'an F, {\it Instantaneous frequency estimation and representation of the audio signal through
      Complex Wavelet Additive Synthesis,} International Journal of Wavelets, Multiresolution and Information Processing, 2014, 12(03):1450030.


\bibitem{LDQ1} H.C. Li, G.T. Deng, T. Qian, {\it Fourier Spectrum Characterizations of $ H^{p} $ Spaces on Tubes Over Cones for $1\leq p\leq \infty,$} Complex Analysis and Operator Theory, (2017). https://doi.org/10.1007/s11785-017-0737-6.

    \bibitem{LDQ2} H. C. Li, G. T.Deng, T. Qian, {\it Hardy space decomposition of on the unit circle: 0<p<1,} Complex Variables and Elliptic Equations: An International Journal, 2016, 61(4): 510-523.


\bibitem{LFZ} Y. Lei, Y. Fang, and L.M. Zhang. {\it Iterative learning control for discrete linear system with wireless transmission based on adaptive fourier decomposition,} Control Conference (CCC), 2017 36th Chinese IEEE, 2017.

\bibitem{LJC} Y. Liang Y, L.-M. Jia, G. Cai, {\it A new approach to diagnose rolling bearing faults based on AFD,} Proceedings of the
    2013 International Conference on Electrical and Information Technologies for Rail Transportation-Volume II, Springer.


\bibitem{LMcQ} C. Li, A. McIntosh, T. Qian, {\it Clifford algebras, Fourier transforms, and singular Convolution operators on Lipschitz surfaces, } Revista Matematica Iberoamericana,1994, 10(3): 665-695.

    \bibitem{LMcS} C. Li, A. McIntosh, S. Semmes, {\it Convolution Singular Integrals on Lipschitz Surfaces,} Journal of the American Mathematical Society, 1992: 455-481.

\bibitem{LQM1} S. Li, T. Qian, W-X. Mai, {\it Sparse Reconstruction of Hardy Signal And Applications to Time-Frequency Distribution,} accepted to appear in International Journal of Wavelets, Multiresolution and Information Processing.

    \bibitem{Ly} A. Lyzzaik, {\it On a conjecture of M.S. Robertson,} Proc. Am. Math. Soc. 91, 108��?10 (1984).

    \bibitem{Na} M. Nahon, {\it Phase Evaluation and Segmentation,} Ph.D. Thesis, Yale University, 2000.

        \bibitem{McQ1} A. McIntosh, T. Qian, {\it Convolution singular integrals on Lipschitz curves,} Springer-Verlag, Lecture Notes in Maths 1494 (1991) 142--162.

            \bibitem{McQ2} A. McIntosh, T. Qian, {\it Lp Fourier multipliers along Lipschitz curves,}  Transactions of The American Mathematical Society, 1992, 333(1): 157-176.

 \bibitem{MF} J. Mashreghi, E. Fricain, {\it Blaschke products and their applications,} Springer, 2013.

    \bibitem{MaiQS1} W.-X. Mai, T. Qian and S. Saitoh, {\it Adaptive Decomposition of Functions with Reproducing Kernels,} in preparation.

\bibitem{MQ1} W. Mi and T. Qian, {\it Frequency Domain Identification: An Algorithm Based On Adaptive Rational Orthogonal System, Automatica,} 48(6). pp. 1154-1162.

\bibitem{MQM} Y. Mo, T. Qian, W. Mi, {\it Sparse Representation in Szego Kernels through Reproducing Kernel Hilbert Space Theory with Applications,} International Journal of Wavelet, Multiresolution and Information Processing, 2015, 13(4), 1550030, 20pp.


\bibitem{MQW}  W. Mi, T. Qian, F. Wan, A Fast Adaptive Model Reduction
Method Based on Takenaka-Malmquist Systems, by W. Mi, T. Qian and
F. Wan, Systems and Control Letters, Volume 61, Issue 1, January
2012, Page 223-230.

\bibitem{MS} F.E. Mozes, J. Szalai, {\it Computing the instantaneous frequency for an ECG signal,} Scientific Bulletin of the $\lq\lq$Petru Maior" University of Targu Mures, 2012, 9(2):28.

\bibitem{Pe} A. Perotti A, {\it Directional quaternionic Hilbert operators,} Hypercomplex analysis, Birkh\"user Basel, 2008:235-258.

\bibitem{Pi} B. Picinbono, {\it On instantaneous amplitude and phase of signals,}
 IEEE Transactions on Signal Processing 1997; 45(3):552--560.

\bibitem{Q11} T. Qian, {\it Singular integrals with holomorphic kernels and Fourier multipliers on star-shape Lipschitz curves,} Studia Mathematica, 1997, 123(3): 195-216.

 \bibitem{Q17} T. Qian, {\it Characterization of boundary values of functions in Hardy spaces with applications in signal analysis,} Journal of Integral Equations and Applications, 2005, 17(2): 159-198.

 \bibitem{Q18} T. Qian, {\it Analytic Signals and Harmonic Measures,} Journal of Mathematical Analysis and Applications, 2006, 314(2): 526-536.


 \bibitem{Q19} T. Qian, {\it Mono-components for decomposition of signals,}  Mathematical Methods in the Applied Sciences, 2006, 29(10): 1187-1198.


\bibitem{Q20} T. Qian, {\it Boundary Derivatives of the Phases of Inner and Outer Functions and Applications,} Mathematical Methods in the Applied Sciences, 2009, 32: 253-263.


\bibitem{Q} T. Qian,  {\it Intrinsic mono-component decomposition of functions:
An advance of Fourier theory,} Mathematical Methods in Applied Sciences, 2010, 33, 880-891, DOI: 10.1002/mma.1214.

\bibitem{Q2D}  T. Qian, {\it Two-Dimensional Adaptive Fourier Decomposition, } Mathematical Methods in the Applied Sciences, 2016, 39(10) : 2431-2448.

    \bibitem{Qbook}  T. Qian, {\it Adaptive Fourier Decomposition: A Mathematical Method Through Complex Analysis, Harmonic Analysis and Signal Analysis,} the Chinese Science Press (in Chinese), 2015.

        \bibitem{Qbook2} T. Qian and P.-T. Li, {\it Singular Integrals and Fourier Theory,} the Chinese Science Press (in Chinese), 2017.



\bibitem{Q15} T. Qian, {\it Fourier analysis on starlike Lipschitz surfaces, } Journal of Functional Analysis, 2001, 183: 370-412.


\bibitem{Qcyclic} T. Qian, {\it  Cyclic AFD Algorithm for Best Approximation by Rational Functions of Given Order,} accepted to appear in Mathematical Methods in the Applied Sciences.

    \bibitem{QChen} T. Qian and Q-H. Chen, {Rational Orthogonal Systems are Schauder Bases,} accepted by Complex Variables and Elliptic Equations.

    \bibitem{QCL} T. Qian, Q.-H. Chen, L.-Q. Li, {\it Analytic unit quadrature signals with non-linear phase,} Physica D: Nonlinear Phenomena, 303 (2005), 80-87.

        \bibitem{QCT} T. Qian, Q. H. Chen, L.H. Tan, {\it Rational Orthogonal Systems are Schauder Bases,} Complex Variables and Elliptic Equations, 2014, 59(6): 841-846.


        \bibitem{QHuang} T. Qian and J.-S. Huang, {\it AFD on the $n$-Torus,} in preparation.

\bibitem{QHLW} T. Qian, I. T. Ho, I. T. Leong, Y. B. Wang, {\it Adaptive decomposition of functions into pieces of non-negative instantaneous frequencies,} International Journal of Wavelets, Multiresolution and Information Processing, 8 (2010), no. 5, 813-833.

    \bibitem{QLS} T. Qian, H. Li, M. Stessin ,{\it  Comparison of Adaptive Mono-component Decompositions, } Nonlinear Analysis: Real World Applications, Volume 14, Issue 2, April 2013, Pages 1055~1074.

\bibitem{QT1} T. Qian, L.-H. Tan, {\it Characterizations of Mono-components: the Blaschke and Starlike types, } Complex Analysis and Operator Theory, 2015, 1-17, DOI 10.1007/s11785-015-0491-6.

        \bibitem{QT2} T. Qian, L. H. Tan, {\it Backward shift invariant subspaces with applications to band preserving and phase retrieval problems,} Mathematical Methods in the Applied Sciences,2016, 39(6): 1591-1598.



\bibitem{QWa}  T. Qian, Yanbo Wang, \textit{Adaptive Fourier Series-A Variation of Greedy
Algorithm}, Advances in Computational Mathematics, 34(2011), no.3,
279-293.

\bibitem{QWe} T. Qian and E. Wegert,
{\it Optimal Approximation by Blaschke Forms,} Complex Variables
and Elliptic Equations, Volume 58, Issue 1, 2013, page 123-133.

\bibitem{QSW} T. Qian, W. Sproessig, J. X. Wang, Adaptive Fourier decomposition of functions in quaternionic Hardy spaces, Mathematical Methods in the Applied Sciences, 2012, 35(1): 43�V64.


\bibitem{QTW} T. Qian, L.H. Tan and Y.B. Wang, {\it Adaptive Decomposition by Weighted Inner Functions: A Generalization of Fourier Serie,} J. Fourier Anal. Appl. 17 (2011), no. 2, 175 ~190.

    \bibitem{QWj} T. Qian and J.-X. Wang, {\it  Adaptive Decomposition of Functions by
Higher Order Szeg\"{o} Kernels I: A Method for Mono-component
Decomposition,} submitted to Acta Applicanda Mathematicae.

\bibitem{QWjz} T. Qian and Jian-Zhong Wang, {\it Gradient Descent Method fpr Best Blaschke-Form Approximation of Function in Hardy Space,} http://arxiv.org/abs/1803.08422.

\bibitem{QSW} T. Qian, W. Sproessig and J.-X. Wang, {\it Adaptive Fourier decomposition of functions in quaternionic Hardy spaces,} Mathematical Methods in the Applied Sciences. (35) 2012, 43 每~64.DOI: 10.1002/mma.1532.


\bibitem{QWY} T. Qian, J. X. Wang, Y. Yang, {\it Matching Pursuits among Shifted Cauchy Kernels in Higher-Dimensional Spaces,} Acta Mathematica Scientia, 2014, 34(3): 660-672.


\bibitem{QWYZ} T. Qian, R. Wang, Y.-S. Xu and H.-Z. Zhang, {\it Orthonormal
Bases with Nonlinear Phase,} Advances in Computational Mathematics 33
(2010), 75-95.

\bibitem{QXYYY} T. Qian, Y. S. Xu, D. Y. Yan, L. X. Yan, B. Yu, Fourier Spectrum Characterization of Hardy Spaces and Applications, Proceedings of the American Mathematical Society, 2009, 137(3): 971-980.

 \bibitem{QY} T. Qian and Y. Yang,  {\it  Hilbert Transforms on the Sphere With the Clifford Algebra Setting,} Journal of Fourier Analysis and Applications,  (2009) 15: 753-774. DOI: 10.1007/s00041-009-9062-4.


\bibitem{QZL} T. Qian, Liming Zhang and Zhi-Xiong Li, {\it Algorithm of Adaptive Fourier Decomposition,}
IEEE Transaction on Signal Processing, Dec., 2011, Volume 59,
Issue 12, page 5899-5902.

\bibitem{QD} W. Qu and P. Dang, {\it Rational Approximation in the Bergman Spaces,} http://arxiv.org/abs/1803.04609.

   \bibitem{SA} L. Salomon, {\it Analyse de l'anisotropie dans des images textur\'ees} 2016.

    \bibitem{SH1} F. Sakaguchi, M. Hayashi, {\it General theory for integer-type algorithm for higher order differential equations,} Numerical Functional Analysis and Optimization, 2011, 32(5):541-582.

\bibitem{SH2} F. Sakaguchi, M. Hayashi, {\it Differentiability of eigenfunctions of the closures of differential operators with rational coefficient functions,} arXiv:0903.4852(2009).

\bibitem{SH3} F. Sakaguchi, M. Hayashi, {\it Practical implementation and error bound of integer-type algorithm for higher-order differential
       equations,}  Numerical Functional Analysis and Optimization, 2011, 32(12): 1316-1364.

\bibitem{SH4} F. Sakaguchi, M. Hayashi, {\it Integer-type algorithm for higher order differential equations by smooth wavepackets,} arXiv:0903.4848(2009).


\bibitem{SQSW} D. Schepper, T. Qian, F. Sommen, J. X. Wang, {\it Holomorphic Approximation of $L_2$-functions on the Unit Sphere in R3,} Journal of  Mathematical Analysis and Applications, 2014, 416(2), 659-671.

\bibitem{SV} RC. Sharpley and V. Vatchev, {\it Analysis of intrinsic mode functions,}  Constructive Approximation 2006; 24:17--47.

    \bibitem{SW} E.M. Stein and G. Weiss, {\it Introduction to Fourirer Analysis on Euclidean Spaces, } Princeton University Press, Princeton, New Jersey, 1971.

    \bibitem{Tan-Shen-Yang}
L.H. Tan, L.X.Shen and L.-H. Yang,
{\it Rational orthogonal bases satisfying the Bedrosian Identity,} Advances in Computational Mathematics (2010)33:285-303.

\bibitem{TYH} L.H. Tan, L.H. Yang, D.R. Huang, {\it The structure of instantaneous frequencies of periodic analytic
signals,} Sci. China Math. 53(2), 347��?55 (2010).

   \bibitem{TQ1} L.-H. Tan, T. Qian, {\it Backward Shift Invariant Subspaces With Applications to Band
Preserving and Phase Retrieval Problems,}

 \bibitem{TQ0} L.H. Tan, T. Qian, {\it Extracting Outer Function Part from Hardy Space Function,} Science China Mathematics, 2017, 60 (11): 2321-2336.

\bibitem{TQC} L. H. Tan, T. Qian, Q. H. Chen, {\it New aspects of Beurling�VLax shift invariant subspaces,} Applied Mathematics and Computation, 2015, 256, 257-266.

\bibitem{VVV} V. Vatchev, {\it A class of intrinsic trigonometric mode polynomials,} International Conference Approximation Theory, Springer, Cham, 2016:361-373.



\bibitem{Daniel} D. V. Vliet, {\it Analytic signals with non-negative instantaneous frequency,} Journal of Integral Equations and Applications, Vol 21, No.1, Spring 2009, 95-111.

\bibitem{Wa} J.L. Walsh, {\it Interpolation and Approximation by Rational Functions in the Complex Plane,} American Mathematical Society: Providence, RI, 1969.

\bibitem{WangSL} S.L. Wang, {\it Simple Proofs of the Bedrosian Equality for the Hilbert Transform,} Science in China, Series A: Mathematics, 2009, 52(3): 507-510.

    \bibitem{WQ1} J. X. Wang, T. Qian, {Approximation of monogenic functions by higher order Szeg\"o kernels on the unit ball and the upper half space,} Sciences in China: Mathematics, 2014, 57(9), 1785-1797.

\bibitem{WW} G. Weiss and M. Weiss, {\it A derivation of the main results of the theory of Hp-spaces,} Rev. Un. Mat. Argentina 20, 1962, 63-71.


    \bibitem{WW1} W. Wu, {\it Applications in Digital Image Processing of Octonions Analysis and the Qian Method,} South China Normal University, 2014.

    \bibitem{WCW},  Z.  Wang, J. N. da Cruz, and F. Wan. {\it Adaptive Fourier decomposition approach for lung-heart sound separation,} Computational Intelligence and Virtual Environments for Measurement Systems and Applications (CIVEMSA), 2015 IEEE International Conference on. IEEE, 2015.

    \bibitem{W1} M.Z. Wu, Y. Wang, X.-M. Li, {\it Fast Algorithm of The Qian Method in Digital Watermarking,} [J]. Computer Engineering and Desining, 2016.

        \bibitem{W2}  M.Z. Wu, Y. Wang, X.-M. Li, {\it Improvement of 2D Qian Method and its Application in Image Denoising,} South China Normal University, 2016.

    \bibitem{Xu} Y.S. Xu,  {\it Private comminication}, 2005.

    \bibitem{YQS2} Y. Yang, T.  Qian, F. Sommen, {\it Phase Derivative of Monogenic Signals in Higher Dimensional Spaces,} Complex Analysis and Operator Theory, 2012, 6(5), 987-1010.



   \bibitem{YuZ} B. Yu and H.Z. Zhang, {\it The Bedrosian Identity and Homogeneous
Semi-convolution Equations,} Journal of Integral Equations and Applications 20 (2008), 527-568.

\bibitem{Z2} L. Zhang, {\it A New Time-Frequency Speech Analysis Approach Based On Adaptive Fourier Decomposition,} World Academy of
Science, Engineering and Technology, International Journal of Electrical, Computer, Energetic, Electronic and Communication
Engineering, 2013.

\bibitem{Z1} L.-M. Zhang, N. Liu, P. Yu, {\it A novel instantaneous frequency algorithm and its application in stock index movement prediction,}
     IEEE Journal of Selected Topics in Signal Processing, 2012, 6(4): 311-318.

\bibitem{ZQDM} L. M. Zhang, T. Qian, W. X. Mai, P. Dang, {\it Adaptive Fourier decomposition-based Dirac type time-frequency distribution,} Mathematical Methods in the Applied Sciences, 2017, 40(8), 2815-2833.




\end{thebibliography}
\end{document}